\documentclass[12pt]{amsart}
\usepackage{amsfonts,amsmath,amssymb}
\newtheorem{Lemma}{Lemma}[section]
\newtheorem{Theorem}[Lemma]{Theorem}
\newtheorem{Prop}[Lemma]{Proposition}

\newtheorem{Cor}[Lemma]{Corollary}

\newcommand{\pf}{\medskip\noindent{\sc Proof: }}
\newcommand{\epf}{$\Box$}

\newcommand{\PH}{\!\!\!\!\!\!\!\!\!\!\!\!\!\! \phantom{\begin{array}{r} a /a \end{array}}}
\newcommand{\UO}{\underline{o}}
\newcommand{\UOP}{\underline{op}}

\newcommand\id{{\operatorname{id}}}

\newcommand\YDG{^{\Gamma}_{\Gamma}\mathcal{YD}}
\newcommand\YDL{^{\Lambda}_{\Lambda}\mathcal{YD}}

\newcommand\G{\Gamma}

\title[Biproducts and Two-Cocycle Twists]{Biproducts and Two-Cocycle Twists of Hopf Algebras}
\author{David E. Radford}
\thanks{Research by the first author partially supported by NSA Grant H98230-04-1-0061. A
significant amount of work on this paper and its sequel was done during his visits to
the Mathematisches Institut der Ludwig-Maximilians-Universit\"{a}t M\"{u}nchen during
June of 2003 and May of 2004 and during the visits of the second author to the
Department of Mathematics, and Statistics, and Computer Science at the University of
Illinois at Chicago during September 2003 and March 2005. The first author expresses
his gratitude for the hospitality and support he received from the Institut and the
second expresses for the same he received from UIC}
\address{University of Illinois at Chicago \\
Department of Mathematics, Statistics and \\ Computer Science (m/c 240) \\
801 South Morgan Street \\
Chicago, IL   60608-7045} \email{radford@uic.edu}
\author{Hans J\"{u}rgen Schneider}
\address{Mathematisches Institut \\
Ludwig-Maximilians-Universit\"{a}t M\"{u}nchen \\
Theresienstr. 39 \\
D-80333 M\"{u}nchen, Germany}
\email{Hans-Juergen.Schneider@mathematik.uni-muenchen.de}

\begin{document}
\numberwithin{equation}{section} \maketitle
\date{}
%
%
\begin{abstract}
{ \small  \rm Let $H$ be a Hopf algebra with bijective antipode over a field $k$ and
suppose that $R{\#}H$ is a bi-product. Then $R$ is a bialgebra in the
Yetter--Drinfel'd category ${}_H^H{\mathcal YD}$. We describe the bialgebras
$(R{\#}H)^{op}$ and $(R{\#}H)^o$ explicitly as bi-products $R^{\UOP}{\#}H^{op}$ and
$R^{\UO}{\#}H^o$ respectively where $R^{\UOP}$ is a bialgebra in
${}^{H^{op}}_{H^{op}}{\mathcal YD}$ and  $R^{\UO}$ is a bialgebra in
${}^{H^o}_{H^o}{\mathcal YD}$. We use our results to describe two-cocycle twist
bialgebra structures on the tensor product of bi-products.}
\end{abstract}
\section*{Introduction}\label{SecIntro}
In \cite{DERHJS21} the irreducible representations of a certain class of pointed Hopf
algebras over a field $k$ is parameterized by pairs of characters, or by characters.
These Hopf algebras are two-cocycle twists $H = (U{\otimes}A)^\sigma$ of the tensor
product of pointed Hopf algebras or quotients of them. The twist structures are in
one-one correspondence with bialgebra maps $U \longrightarrow {A^{op}}^{\;o}$. In many
cases the pointed Hopf algebras $U$ and $A$ are bi-products. We are thus led to
consider the multiplicative opposite $(R{\#}H)^{op}$ and the dual $(R{\#}H)^o$ of a
bi-product $R{\#}H$. Generally the multiplicative opposite and dual of a bi-product is
a bi-product. One purpose of this paper is to characterize $(R{\#}H)^{op}$ and
$(R{\#}H)^o$ as bi-products when $H$ has bijective antipode. Recall that pointed Hopf
algebras have bijective antipodes.

Let $H$ be a Hopf algebra over $k$ with bijective antipode and suppose that $R{\#}H$
is a bi-product. Then $R$ is a bialgebra in the category of Yetter--Drinfel'd modules
${}_H^H{\mathcal YD}$. We construct a bialgebra $R^{\UOP}$ in the Yetter--Drinfel'd
category ${}_{H^{op}}^{H^{op}}{\mathcal YD}$ such that $(R{\#}H)^{op} \simeq
R^{\UOP}{\#}H^{op}$. Likewise we construct a bialgebra $R^{\UO}$ in the
Yetter--Drinfel'd category ${}_{H^o}^{H^o}{\mathcal YD}$ such that $(R{\#}H)^o \simeq
R^{\UO}{\#}H^o$. These bialgebra constructions are based on more general procedures:
the construction of an algebra (respectively a coalgebra) $A^{\UOP}$ in
${}_{H^{op}}^{H^{op}}{\mathcal YD}$ from an algebra (respectively a coalgebra) $A$ in
${}_H^H{\mathcal YD}$ and the construction of an algebra (respectively a coalgebra)
$A^{\UO}$ in ${}_{H^o}^{H^o}{\mathcal YD}$ from an algebra (respectively a coalgebra)
$A$ in ${}_H^H{\mathcal YD}$.

Important to us is the case $U = {\mathfrak B}(W){\#}k[\G]$ and $A = {\mathfrak
B}(V){\#}k[\Lambda]$, where $\Gamma$ and $\Lambda$ are abelian groups and $({\mathfrak B}(W)$, ${\mathfrak B}(V)$ are Nichols algebras in the categories
${}_{k[\G]}^{k[\G]}{\mathcal YD}$, ${}_{k[\Lambda]}^{k[\Lambda]}{\mathcal YD}$
respectively. These are of primary interest in \cite{DERHJS21}. For these Hopf algebras
an extensive class of bialgebra maps $U \longrightarrow A^{op \, o}$ can be given in
terms of two linear forms $\tau : k[\Lambda]{\otimes}k[\G] \longrightarrow k$ and
$\beta : V{\otimes}W \longrightarrow k$. The forms can easily be produced and thus, in
particular, our results provide a way of generating a large number of two-cocycle
twist bialgebras, that is bialgebra maps $\mathfrak{B}(W){\#}k[\Lambda] \longrightarrow (\mathfrak{B}(V){\#}k[\G])^{op \, o}$ without checking the relations of $\mathfrak{B}(W)$ or $\mathfrak{B}(V)$ which are unknown in general.

For finite abelian groups $\Gamma, V \in {}_{k[\G]}^{k[\G]}{\mathcal YD}$ and finite-dimensional Nichols algebras $\mathfrak{B}(V)$, the dual of $\mathfrak{B}(V) {\#} k[\G]$ was already computed in \cite[Theorem 2.2]{B}. The two-cocycles in Corollary \ref{groupcase} were determined in \cite{AS} for finite-dimensional Nichols algebras with known relations by explicitly checking the relations. However, in \cite{AS} the more general case when $U$ is not coradically graded, that is of the form $\mathfrak{B}(W){\#}k[\Lambda],$ was considered.

This paper is organized as follows. In Section \ref{SecPrelim} we deal with the
somewhat extensive prerequisites for the paper. First we discuss notations for
algebras, coalgebras, and their representations, and then review algebraic objects in
the Yetter--Drinfel'd category ${}_H^H{\mathcal YD}$ of a Hopf algebra $H$ with
bijective antipode in detail for the reader's convenience. This discussion is
important for Sections \ref{SecHop} and \ref{SecHo} where we describe algebra,
coalgebra, and bialgebra constructions in the Yetter--Drinfel'd categories and
${}_{H^{op}}^{H^{op}}{\mathcal YD}$ and ${}_{H^o}^{H^o}{\mathcal YD}$ based on the
their counterparts in ${}_H^H{\mathcal YD}$. These constructions are basic ingredients
in the realization of the multiplicative opposite and dual of a bi-product as a
bi-product.

Certain bilinear forms on objects in Yetter--Drinfel'd categories which play a role in
our construction of our two-cocycle twist Hopf algebras are introduced and studied in
Section \ref{SecFormYD}. In Section \ref{SectionBiProd} we consider morphisms of
bi-products. The multiplicative opposite of a bi-product is characterized as a
bi-product in Section \ref{SubAsmashHOpp} and the dual of a bi-product is
characterized as a bi-product in Section \ref{SubAsmashHO}.

We apply the main results of Sections \ref{SubAsmashHOpp} and \ref{SubAsmashHO} to
describe certain two-cocycle twists on the tensor product $(T{\#}K){\otimes}(R{\#}H)$
of bi-products in Section \ref{SecUABiProducts}. Here $K$ and $H$ are Hopf algebras
with bijective antipodes over $k$. We focus on the basic case when $T = {\mathfrak
B}(W)$ and $R = {\mathfrak B}(V)$ are Nichols algebras on finite-dimensional
Yetter--Drinfel'd modules. In the last section we consider the our basic case when $K
= k[\Lambda]$ and $H = k[\G]$ are group algebras of abelian groups. Results here fit
nicely into the discussion of \cite{DERHJS21}.

We denote the antipode of a Hopf algebra $H$ over $k$ by $S$. Any one of
\cite{Abe,LamRad,Mont,SweedlerBook} will serve as a Hopf algebra reference for this
paper.  Throughout $k$ is a field and all vector spaces are over $k$. For vector
spaces $U$ and $V$ we will drop the subscript $k$ from ${\rm End}_k(V)$, ${\rm
Hom}_k(U, V)$, and $U{\otimes}_kV$. We denote the identity map of $V$ by $I_V$. For a
non-empty subset $S$ of the dual space $V^*$ we let $S^\perp$ denote the subspace of
$V$ consisting of the common zeros of the functionals in $S$. For $p \in U^*$ and $u
\in U$ we denote the evaluation of $p$ on $u$ by $p(u)$ or ${<}p, u{>}$.
\section{Preliminaries}\label{SecPrelim}
A good deal of prerequisite material is needed for this paper. We discuss general
notation first then review specific topics in detail.

For a group $G$ we let $\widehat{G}$ denote the group of characters of $G$ with values
in $k$. $H = k[G]$ denotes the group algebra of $G$ over $k$ which is Hopf algebra
arising in most applications in this paper. For a Hopf algebra $H$ over $k$ we denote
the group of grouplike elements of $H$ by $G(H)$ as usual.

Let $(A, m, \eta)$ be an algebra over $k$, which we shall usually denote by $A$.
Generally we represent algebraic objects defined on a vector space by their underlying
vector space. Observe that $(A, m^{op}, \eta)$ is an algebra over $k$, where $m^{op} =
m{\circ}\tau_{A, A}$. We denote $A$ with this algebra structure by $A^{op}$ and write
$a^{op} = a$ for elements of $A^{op}$. Thus $a^{op}b^{op} = (ba)^{op}$ for all $a, b
\in A$. This notation is very useful for computations in Yetter--Drinfel'd categories
discussed below involving certain algebra constructions. We denote the category of
left (respectively right) $A$-modules and module maps by ${}_A{\mathcal M}$
(respectively ${\mathcal M}_A$). If ${\mathcal C}$ is a category, by abuse of notation
we will write $C \in {\mathcal C}$ to indicate that $C$ is an object of ${\mathcal
C}$.

Let $M$ be a left $A$-module. Then $M^*$ is a right $A$-module under the transpose
action  which is given by $(m^*{\cdot}a)(m) = m^*(a{\cdot}m)$ for all $m^* \in M^*$,
$a \in A$, and $m \in M$. Likewise if $M$ is a right $A$-module then $M^*$ is a left
$A$-module where $(a{\cdot}m^*)(m) = m^*(m{\cdot}a)$ for all $a \in A$, $m^* \in M^*$,
and $m \in M$.

Suppose $B$ is an algebra also over $k$, let $N$ be a left $B$-module, and suppose
that $\varphi : A \longrightarrow B$ is an algebra map. Then a linear map $f : M
\longrightarrow N$ is $\varphi$-linear if $f(a{\cdot}m) = \varphi(a){\cdot}f(m)$ for
all $a \in A$ and $m \in M$. There is a way of expressing the last equation in terms
of $A$-module maps. Note that $N$ is a left $A$-module by pullback along $\varphi$.
Thus $f$ is $\varphi$-linear if and only if $f$ is a map of left $A$-modules.

Let $(C, \Delta, \epsilon)$ be a coalgebra over $k$, which we usually denote by $C$.
At times it is convenient to denote the coproduct $\Delta$ by $\Delta_C$. Generally we
use a variant on the Heyneman-Sweedler notation for the coproduct and write $\Delta
(c) = c_{(1)}{\otimes}c_{(2)}$ to denote $\Delta(c) \in C{\otimes}C$ for $c \in C$.
Note that $(C, \Delta^{cop}, \epsilon)$ is a coalgebra over $k$, where $\Delta^{cop} =
\tau_{C, C}{\circ}\Delta$. We let $C^{cop}$ denote the vector space $C$ with this
coalgebras structure and sometimes write $c^{cop} = c$ for elements of $C^{cop}$. With
this notation $c^{cop}_{\;\;\;(1)}{\otimes}c^{cop}_{\;\;\;(2)} =
c_{(2)}{\otimes}c_{(1)}$ for all $c \in C^{cop}$.

Suppose that $(M, \delta)$ is a left $C$-comodule. For $m \in M$ we use the notation
$\delta (m) = m_{(-1)}{\otimes}m_{(0)}$ to denote $\delta (m) \in C{\otimes}M$. If
$(M, \delta)$ is a right $C$-comodule we denote $\delta (m) \in M{\otimes}C$ by
$\delta (m) = m_{(0)}{\otimes}m_{(1)}$.  Observe that our coproduct and comodule
notations do not conflict.

We make an exception to our coproduct notation described above for coalgebras in
Yetter--Drinfel'd categories, in which case we write $\Delta(c) =
c^{(1)}{\otimes}c^{(2)}$ for $c \in C$. See Section \ref{SecYDCatBiProd}.

Suppose that $M$ is a left $C$-comodule, $D$ is a coalgebra over $k$, $N$ is a left
$D$-comodule, and $\varphi : C \longrightarrow D$ is a coalgebra map. Then a linear
map $f : M \longrightarrow N$ is left $\varphi$-colinear if $\varphi
(m_{(-1)}){\otimes}f(m_{(0)}) = f(m)_{(-1)}{\otimes}f(m)_{(0)}$ for all $m \in M$.
There is a way of expressing the last equation in terms of $D$-comodule maps. Note
that $M$ is a left $D$-comodule by pushout along $\varphi$. Thus $f$ is
$\varphi$-colinear if and only if $f$ is a map of left $D$-modules. We use the
terminology $\varphi$-linear and colinear as shorthand for $\varphi$-linear and
$\varphi$-colinear.

Bilinear forms play an important role in this paper. We will think of them in terms of
linear forms $\beta : U{\otimes}V \longrightarrow k$ and will often write $\beta (u,
v)$ for $\beta (u{\otimes}v)$. Note that $\beta$ determines linear maps $\beta_\ell :
U \longrightarrow V^*$ and $\beta_r : V \longrightarrow U^*$ where $\beta_\ell (u)(v)
= \beta (u, v) = \beta_r(v)(u)$ for all $u \in U$ and $v \in V$. The form $\beta$ is
left (respectively right) non-singular if $\beta_\ell$ (respectively $\beta_r$) is
one-one and $\beta$ is non-singular if it is both left and right non-singular.

For subspaces $X \subseteq U$ and $Y\subseteq V$ we define subspaces $X^\perp
\subseteq V$ and $Y^\perp \subseteq U$ by
$$
X^\perp = \{ v \in V \,|\, \beta (X, v) = (0)\,\} \quad \mbox{and} \quad Y^\perp = \{
u \in U \,|\, \beta (u, Y) = (0)\,\}.
$$
Note that there is a form $\overline{\beta} : U/V^\perp{\otimes}V/U^\perp
\longrightarrow k$ uniquely determined by
$\overline{\beta}{\circ}(\pi_{V^\perp}{\otimes}\pi_{U^\perp}) = \beta$, where
$\pi_{V^\perp} : U \longrightarrow U/V^\perp$ and $\pi_{U^\perp} : U \longrightarrow
V/U^\perp$ are the projections. Observe that $V^\perp = {\rm Ker}\,\beta_\ell$,
$U^\perp = {\rm Ker}\,\beta_r$, and that $\overline{\beta}$ is non-singular.
\subsection{Two-Cocycle Twist Bialgebras}\label{SubSecTwist}
Let $A$ be a bialgebra over $k$. A two-cocycle for $A$ is a convolution invertible
linear form $\sigma : A{\otimes}A \longrightarrow k$ which satisfies
$$
\sigma (x_{(1)}, y_{(1)})\sigma(x_{(2)}y_{(2)}, z) = \sigma(y_{(1)}, z_{(1)})\sigma(x,
y_{(2)}z_{(2)})
$$
for all $x, y, z \in A$. If $\sigma$ is a two-cocycle for $A$ then $A^\sigma$ is a
bialgebra, where $A^\sigma = A$ as a coalgebra and multiplication $m^\sigma :
A{\otimes}A \longrightarrow A$ is given by
$$
m^\sigma (x{\otimes}y) = \sigma(x_{1)}, y_{(1)})x_{(2)}y_{(2)}\sigma^{-1}(x_{(3)},
y_{(3)})
$$
for all $x, y \in A$.

Let $U$ and $A$ be bialgebras over $k$ and suppose that $\tau : U{\otimes}A
\longrightarrow k$ is a linear form. Consider the axioms:
\begin{enumerate}
\item[{\rm (A.1)}] $\tau (u, aa') = \tau (u_{(2)}, a)\tau
(u_{(1)}, a')$ for all $u \in U$ and $a, a' \in A$; \item[{\rm (A.2)}] $\tau (1, a) =
\epsilon(a)$ for all $a \in A$; \item[{\rm (A.3)}] $\tau (uu', a) = \tau (u,
a_{(1)})\tau (u', a_{(2)})$ for all $u, u'  \in U$ and $a \in A$; \item[{\rm (A.4)}]
$\tau (u, 1) = \epsilon(u)$ for all $u \in U$.
\end{enumerate}
Axioms (A.1)--(A.4) are equivalent to
\begin{equation}\label{EqTauUAopo}
\mbox{$\tau_\ell (U) \subseteq A^o$ and $\tau_\ell : U \longrightarrow A^{o\, cop} =
A^{op \, o}$ is a bialgebra map}
\end{equation}
and they are also equivalent to
\begin{equation}\label{EqTauUAUoop} \mbox{$\tau_r (A) \subseteq U^o$ and $\tau_r : A
\longrightarrow U^{o\, op}$ is a bialgebra map.}
\end{equation}

Suppose that (A.1)--(A.4) hold, $\tau$ is convolution invertible, and define a linear
form $\sigma : (U{\otimes}A){\otimes}(U{\otimes}A) \longrightarrow k$ by $\sigma
(u{\otimes}a, u'{\otimes}a') = \epsilon (a)\tau(u', a)\epsilon(a')$ for all $u, u' \in
U$ and $a, a' \in A$. Then $\sigma$ is a two-cocycle. We denote the two-cocycle twist
bialgebra structure on the tensor product bialgebra $U{\otimes}A$ by $H =
(U{\otimes}A)^\sigma$. Observe that
$$
(u{\otimes}a)(u'{\otimes}a') = u\tau(u'_{(1)},
a_{(1)})u'_{(2)}{\otimes}a_{(2)}\tau^{-1}(u'_{(3)}, a_{(3)})a'
$$
for all $u, u' \in U$ and $a, a' \in A$.

Suppose that (A.1)--(A.4) hold for the linear form $\tau : U{\otimes}A \longrightarrow
k$. Then $\tau$ is invertible if $U$ has an antipode $S$ or if $A^{op}$ has an
antipode $\varsigma$. In the first case $\tau^{-1}(u, a) = \tau(S(u), a)$, and in the
second $\tau^{-1}(u, a) = \tau(u, \varsigma(a))$, for all $u \in U$ and $a \in A$. See
\cite[Lemma 1.2]{DERHJS21}.

As noted in \cite[Section 1]{DERHJS21}, the quantum double provides an important
example of a two-cocycle twist bialgebra. This example is described in \cite{DoiTak}
where two-cocycle twist bialgebras are defined and discussed.

Let $U$, $\overline{U}$ and $A$, $\overline{A}$ be algebras over $k$. Suppose further
that $\tau :U{\otimes}A \longrightarrow k$ and $\overline{\tau} :
\overline{U}{\otimes}\overline{A} \longrightarrow k$ are convolution invertible linear
forms satisfying (A.1)--(A.4). Set $H = (U{\otimes}A)^{\sigma}$ and $\overline{H} =
(\overline{U}{\otimes}\overline{A})^{\overline{\sigma}}$. Suppose that $f : U
\longrightarrow \overline{U}$ and $g : A \longrightarrow \overline{A}$ are bialgebra
maps such that $\overline{\tau}(f(u), g(a)) = \tau (u, a)$ for all $u \in U$ and $a
\in A$. Then $f{\otimes}g : H \longrightarrow \overline{H}$ is a bialgebra map.

\subsection{Yetter--Drinfel'd Categories and Their Algebras, Coalgebras,
and Bialgebras}\label{SecYDCatBiProd} Here we organize well-known material for the
reader's convenience and for our use in later sections. See \cite{ASSurvey} in
particular.

Let $H$ be a bialgebra over $k$ and let ${}_H^H{\mathcal YD}$ be the category whose
objects are triples $(M, {\cdot}, \delta)$, where $(M, {\cdot})$ is a left $H$-module,
$(M, \delta)$ is a left $H$-comodule, compatible in the sense that
\begin{equation}\label{EqYDComp}
h_{(1)}m_{(-1)}{\otimes}h_{(2)}{\cdot}m_{(0)} =
(h_{(1)}{\cdot}m)_{(-1)}h_{(2)}{\otimes}(h_{(1)}{\cdot}m)_{(0)}
\end{equation}
for all $h \in H$ and $m \in M$, and whose morphisms $(M, {\cdot}, \delta)
\longrightarrow (M', {\cdot}', \delta')$ are maps $f : M \longrightarrow M'$
simultaneously of left $H$-modules and of left $H$-comodules.  We follow the
convention of referring to an object of ${}_H^H{\mathcal YD}$ as a Yetter--Drinfel'd
module \cite{ASSurvey}. If $H$ has an antipode $S$ then (\ref{EqYDComp}) is equivalent
to
\begin{equation}\label{EqYDCompAntipode}
\delta (h{\cdot}m) = h_{(1)}m_{(-1)}S(h_{(3)}){\otimes}h_{(2)}{\cdot}m_{(0)}
\end{equation}
for all $h \in H$ and $m \in M$,  in practice a very useful formulation of the
compatibility condition (\ref{EqYDComp}).

The category ${}_H^H{\mathcal YD}$ has a monoidal structure, where $k$ is given the
left $H$-module structure $h{\cdot}1_k = \epsilon (h)$ for all $h \in H$ and left
$H$-comodule structure determined by $\delta (1_k) = 1_H{\otimes}1_k$, and the tensor
product of objects $M, N \in {}_H^H{\mathcal YD}$ is $M{\otimes}N$ as a vector space
with left $H$-module structure given by $h{\cdot}(m{\otimes}n) =
h_{(1)}{\cdot}m{\otimes}h_{(2)}{\cdot}n$ and left $H$-comodule structure given by
$\delta (m{\otimes}n) = m_{(-1)}n_{(-1)}{\otimes}(m_{(0)}{\otimes}n_{(0)})$ for all $h
\in H$, $m \in M$, and $n \in N$. When $H$ is a Hopf algebra ${}_H^H{\mathcal YD}$ is
a braided monoidal category with braiding $\sigma_{M, N} : M{\otimes}N \longrightarrow
N{\otimes}M$ for objects $M, N \in {}_H^H{\mathcal YD}$ determined by $\sigma_{M,
N}(m{\otimes}n) = m_{(-1)}{\cdot}n{\otimes}m_{(0)}$ for all $m \in M$ and $n \in N$.

Let $(A, m, \eta)$ be an algebra in ${}_H^H{\mathcal YD}$. Then $(A, m^{op}, \eta)$ is
a well, where $m^{op} = m{\circ}\sigma_{A, A}$. Thus $a^{op}b^{op} =
(a_{(-1)}{\cdot}b)a_{(0)}$ for all $a, b \in A$. We denote the object $A$ with this
algebra structure by $A^{op}$. If $B$ is also an algebra in ${}_H^H{\mathcal YD}$ then
$A{\otimes}B$ is an algebra in ${}_H^H{\mathcal YD}$, where $\eta_{A{\otimes}B} =
\eta_A{\otimes}\eta_B$ and $m_{A{\otimes}B} =
(m_A{\otimes}m_B){\circ}(I_A{\otimes}\sigma_{A, B}{\otimes}I_B)$. We write
$A\underline{\otimes}B$ for $A{\otimes}B$ with this algebra structure and
$a\underline{\otimes}b = a{\otimes}b$ for tensors. By definition
$$
(a\underline{\otimes}b)(a'\underline{\otimes}b') =
a(b_{(-1)}{\cdot}a')\underline{\otimes}b_{(0)}b'
$$
for all $a, a' \in A$ and $b, b' \in B$. Observe that the object $k$ with its usual
$k$-algebra structure is an algebra in ${}_H^H{\mathcal YD}$.

Suppose $(C, \Delta, \epsilon)$ is a coalgebra in ${}_H^H{\mathcal YD}$. We shall
write $\Delta (c) = c^{(1)}{\otimes}c^{(2)}$ for $c \in C$. Observe that $(C,
\Delta^{cop}, \epsilon)$ is a coalgebra in ${}_H^H{\mathcal YD}$, where $\Delta^{cop}
= \sigma_{C, C}{\circ}\Delta$, or equivalently $\Delta^{cop}(c) = c^{(1)}_{\;\;\;
(-1)}{\cdot}c^{(2)}{\otimes}c^{(1)}_{\;\;\; (0)}$ for all $c \in C$. We denote the
object $C$ with this coalgebra structure by $C^{cop}$. If $D$ is a coalgebra in
${}_H^H{\mathcal YD}$ also then $C{\otimes}D$ is a coalgebra in ${}_H^H{\mathcal YD}$,
where $\epsilon_{C{\otimes}D} = \epsilon_C{\otimes}\epsilon_D$ and
$\Delta_{C{\otimes}D} = (I_C{\otimes}\sigma_{C,
D}{\otimes}I_D){\circ}(\Delta_C{\otimes}\Delta_D)$. We write $C\overline{\otimes}D$
for $C{\otimes}D$ with this coalgebra structure. By definition
$$
\Delta(c\overline{\otimes}d) =
(c^{(1)}\overline{\otimes}c^{(2)}_{\;\;\;(-1)}{\cdot}d^{(1)})
{\otimes}(c^{(2)}_{\;\;\;(0)}\overline{\otimes}d^{(2)})
$$
for all $c \in C$ and $d \in D$. Observe that the object $k$ with its usual
$k$-coalgebra structure is a coalgebra in ${}_H^H{\mathcal YD}$.

Let $R \in {}_H^H{\mathcal YD}$ be an algebra and a coalgebra in the category. Then
$\Delta : R \longrightarrow R\underline{\otimes}R$ and $\epsilon : R \longrightarrow
k$ are algebra maps if and only if $m : R\overline{\otimes}R \longrightarrow R$ and
$\eta : k \longrightarrow R$ are coalgebra maps in which case we say that $R$ with its
algebra and coalgebra structure is a bialgebra in ${}_H^H{\mathcal YD}$. If $R$ is a
bialgebra in ${}_H^H{\mathcal YD}$ then $R^{op}$, $R^{cop}$, and therefore $R^{op \;
cop}$, are bialgebras in ${}_H^H{\mathcal YD}$. Observe that the object $k$ with its
usual $k$-algebra and $k$-coalgebra structure is a bialgebra in ${}_H^H{\mathcal YD}$.
If $R$ and $T$ are bialgebras in ${}_H^H{\mathcal YD}$ then the object $R{\otimes}T$
with its algebra structure $R\underline{\otimes}T$ and coalgebra structure
$R\overline{\otimes}T$ is a bialgebra in ${}_H^H{\mathcal YD}$.

The most important bialgebra in ${}_H^H{\mathcal YD}$ for us is the Nichols algebra.
Let $M \in {}_H^H{\mathcal YD}$ and consider the tensor $k$-algebra $T(M) =
k{\oplus}M{\oplus}(M{\otimes}M) {\oplus} \cdots = \bigoplus_{\ell = 0}^\infty
M^{\otimes \, \ell}$ on the vector space $M$. Regard $T(M)$ as an object of
${}_H^H{\mathcal YD}$, where
$$
h{\cdot}(m_1{\otimes} \cdots {\otimes}m_r) = h_{(1)}{\cdot}m_1{\otimes} \cdots
{\otimes}h_{(r)}{\cdot}m_r
$$
for all $h \in H$ and $m_1$, $\ldots$, $m_r \in M$ and
$$
\delta (m_1{\otimes} \cdots {\otimes}m_r) = m_{1 \, (-1)}\cdots m_{r \,
(-1)}{\otimes}(m_{1 \, (0)}{\otimes} \cdots {\otimes}m_{r \, (0)})
$$
for all $m_1$, $\ldots$, $m_r \in M$, describe the left $H$-module and left
$H$-comodule structures respectively. Let $i : M \longrightarrow T(M)$ be the
inclusion. Then the pair $(i, T(M))$ satisfies the obvious analog in ${}_H^H{\mathcal
YD}$ of the universal mapping property of the tensor algebra of a vector space as a
$k$-algebra. Observe that $T(M) = {\bigoplus}_{\ell = 0}^\infty M^{{\otimes}\, \ell}$
is a graded bialgebra (indeed Hopf algebra) in ${}_H^H{\mathcal YD}$.

The algebra $T(M)$ of ${}_H^H{\mathcal YD}$ is a bialgebra in the category just as the
tensor algebra of a vector space is a $k$-bialgebra. The linear maps $d : M
\longrightarrow T(M)\underline{\otimes}T(M)$ and $e : M \longrightarrow k$ defined by
$d(m) = 1\underline{\otimes}m + m\underline{\otimes}1$ and $e(m) = 0$ respectively for
all $m \in M$ lift to algebra morphisms $\Delta : T(M) \longrightarrow
T(M)\underline{\otimes}T(M)$ and $\epsilon : T(M) \longrightarrow k$ uniquely
determined by $\Delta{\circ}i = d$ and $\epsilon{\circ}i = e$. The structure $(T(M),
\Delta, \epsilon)$ is a coalgebra in ${}_H^H{\mathcal YD}$ and $T(M)$ with its algebra
and coalgebra structure is a bialgebra in ${}_H^H{\mathcal YD}$. The pair $(i, T(M))$
satisfies the following universal mapping property: If $A$ is a bialgebra in
${}_H^H{\mathcal YD}$ and $f : M \longrightarrow A$ is a morphism such that ${\rm
Im}\,f \subseteq P(A)$ then there is a bialgebra morphism $F : T(M) \longrightarrow A$
uniquely determined by $F{\circ}i = f$.

The Nichols algebra ${\mathfrak B}(M)$ has a very simple theoretical description.
Among the graded subobjects $J$ of ${\mathfrak B}(M)$ which are coideals and satisfy
$J{\cap}M = (0)$, there is a unique maximal one $I$. It is easy to see that $I$ is an
ideal of ${\mathfrak B}(M)$; thus the subobject $I$ is a bi-ideal. Consequently the
quotient ${\mathfrak B}(M) = T(M)/I$ is a connected graded bialgebra in the category
${}_H^H{\mathcal YD}$. Observe that ${\mathfrak B}(M)$ is generated as an algebra by
${\mathfrak B}(M)(1) = M$, which is also the space of primitive elements of
${\mathfrak B}(M)$. Since $I{\cap}{\mathfrak B}(M) = (0)$ we may think of $M$ as a
subspace of ${\mathfrak B}(M)$. The pair $(M, ({\mathfrak B}(M))$ satisfies the
following universal mapping property:
\begin{Theorem}\label{ThmUMPBM}
Let $H$ be a Hopf algebra and $M \in {}_H^H{\mathcal YD}$. Then:
\begin{enumerate}
\item[{\rm a)}] ${\mathfrak B}(M)$ is a connected graded bialgebra
in ${}_H^H{\mathcal YD}$ and $M = {\mathfrak B}(1)$ is a subobject which generates
${\mathfrak B}(M)$ as an algebra. \item[{\rm b)}] If $A$ is a connected graded
bialgebra in ${}_H^H{\mathcal YD}$ generated as an algebra by $A(1)$ and $f : A(1)
\longrightarrow M$ is a morphism, then there is a bialgebra morphism $F : A
\longrightarrow {\mathfrak B}(M)$ determined by $F|_{A(1)} = f$.
\end{enumerate}
\end{Theorem}

\pf We need only show part b). By the universal mapping property of the bialgebra
$(T(A(1)), i)$ there is a bialgebra morphism $F : T(A(1)) \longrightarrow A$
determined by $F|_{A(1)} = i$. Since $T(A(1))$ is generated by $A(1)$ as an algebra,
$F$ is an onto morphism of graded bialgebras. Let $J = {\rm Ker}\,F$. Then $J$ is a
sub-object of $T(A(1))$ which is a graded bi-ideal of $T(A(1))$ satisfying
$J{\cap}A(1) = (0)$. Using the universal mapping property again we see that the
morphism $f : A(1) \longrightarrow M$ induces a bialgebra morphism $T(f) : T(A(1))
\longrightarrow T(M)$ determined by $T(f)|_{A(1)} = f|_{A(1)}$. Observe that $T(f)(J)$
is a subobject of $T(M)$ which is a graded bi-ideal of $T(M)$ whose intersection with
$M$ is $(0)$. This means $T(f)(J) \subseteq I$, where the latter is defined above. The
composite $A \simeq T(A(1))/J \longrightarrow T(M)/I = {\mathfrak B}(M)$, where the
second map is defined by $x + J \mapsto T(f)(x) + I$, is our desired bialgebra
morphism $F$. \epf
\medskip

We have noted that ${\mathfrak B}(M)(1)$ is the subspace of primitive elements of
${\mathfrak B}(M)$. Thus ${\mathfrak B}(M)$ is a connected graded primitively
generated bialgebra in ${}_H^H{\mathcal YD}$ with subspace of primitive elements
${\mathfrak B}(M)(1)$. These are defining properties.
\begin{Cor}\label{CorNicChar}
Let $H$ be a Hopf algebra over the field $k$ and suppose that $A$ is a connected
graded primitively generated bialgebra in ${}_H^H{\mathcal YD}$ with subspace of
primitive elements $A(1)$. Then there is an isomorphism of bialgebras $A \simeq
{\mathfrak B}(A(1))$ which extends the identity map $I_{A(1)}$.
\end{Cor}

\pf Let $F : A \longrightarrow {\mathfrak B}(A(1))$ be the bialgebra morphism of part
b) of Theorem \ref{ThmUMPBM} which extends $I_{A(1)}$. Since $A(1)$ generates
${\mathfrak B}(A(1))$ the map $F$ is onto. Now ${\rm Ker}\,F{\cap}P(A) = {\rm
Ker}\,F{\cap}A(1) = (0)$. Generally if $C$ is a connected coalgebra and $f : C
\longrightarrow C'$ is a coalgebra map which satisfies ${\rm Ker}\,f{\cap}P(C) = (0)$
then $f$ is one-one \cite[Lemma 11.0.1]{SweedlerBook}. Thus the onto map $F$ is
one-one. \epf
\medskip

We leave the reader with the exercise of verifying the following corollary to the
theorem above.
\begin{Cor}\label{CorNicKH}
Let $K$ and $H$ be Hopf algebras with bijective antipodes over the field $k$ and let
$\varphi : K \longrightarrow H$ be a bialgebra map. Suppose that $W \in
{}_K^K{\mathcal YD}$, $W \in {}_H^H{\mathcal YD}$, and $f : W \longrightarrow V$ is
$\varphi$-linear and colinear. Then:
\begin{enumerate}
\item[{\rm a)}] There is a map of algebras and coalgebras
${\mathfrak B}(f) : {\mathfrak B}(W) \longrightarrow {\mathfrak B}(V)$ determined by
${\mathfrak B}(f)|_W = f$. Furthermore ${\mathfrak B}(f)$ is $\varphi$-linear and
colinear. \item[{\rm b)}] If $f$ is one-one (respectively onto) then ${\mathfrak
B}(f)$ is one-one (respectively onto).
\end{enumerate}
\epf
\end{Cor}

\section{Associated Constructions in ${}_{H^{op}}^{H^{op}}{\mathcal YD}$
}\label{SecHop}
Throughout this section $H$ has bijective antipode $S$. Starting with objects,
algebras, coalgebras, and bialgebras in ${}_H^H{\mathcal YD}$ we construct
counterparts in ${}_{H^{op}}^{H^{op}}{\mathcal YD}$ which are important for the
analysis of bi-products in Section \ref{SecUABiProducts}. First we start with objects.

Let $(M, {\cdot}, \delta) \in {}_H^H{\mathcal YD}$. Then $(M, {\cdot}_{op}, \delta)
\in {}_{H^{op}}^{H^{op}}{\mathcal YD}$, where
\begin{equation}\label{EqMultOp}
h{\cdot}_{op}m = S^{-1}(h){\cdot}m
\end{equation}
for all $h \in H$ and $m \in M$. We denote  $(M, {\cdot}_{op}, \delta)$ by $M^{\UOP}$.
If $N$ is also an object of ${}_H^H{\mathcal YD}$ and $f : M \longrightarrow N$ is a
morphism, then $f^{\UOP} : M^{\UOP} \longrightarrow N^{\UOP}$ is a morphism, where
$f^{\UOP} = f$.

When $M$ has the structure of an algebra, coalgebra, or bialgebra, then $M^{\UOP} \in
{}_{H^{op}}^{H^{op}}{\mathcal YD}$ does as well. If $(A, m, \eta)$ is an algebra in
${}_H^H{\mathcal YD}$ then $(A^{\UOP}, m^{\UOP}, \eta)$ is an algebra in
${}_{H^{op}}^{H^{op}}{\mathcal YD}$, where
\begin{equation}\label{mOP}
m^{\UOP}(a{\otimes}b) = ba
\end{equation}
for all $a, b \in A$. If $\overline{A}$ is also an algebra in ${}_H^H{\mathcal YD}$
and $f : A \longrightarrow \overline{A}$ is an algebra morphism then $f : A^{\UOP}
\longrightarrow \overline{A}^{\UOP}$ is an algebra morphism. If $(C, \Delta,
\epsilon)$ is a coalgebra in ${}_H^H{\mathcal YD}$ then $(C^{\UOP}, \Delta^{\UOP},
\epsilon)$ is a coalgebra in ${}_{H^{op}}^{H^{op}}{\mathcal YD}$, where
\begin{equation}\label{EqDeltaOP}
\Delta^{\UOP}(c) =
c^{(2)}_{\;\;\;\;(-1)}{\cdot}_{op}c^{(1)}{\otimes}c^{(2)}_{\;\;\;\;(0)}
\end{equation}
for all $c \in C$. If $C'$ is also a coalgebra in ${}_H^H{\mathcal YD}$ and $f : C
\longrightarrow \overline{C}$ is a coalgebra morphism then $f : C^{\UOP}
\longrightarrow \overline{C}^{\UOP}$ is a coalgebra morphism. If $(R, m, \eta, \Delta,
\epsilon)$ is a bialgebra in ${}_H^H{\mathcal YD}$ then $(R^{\UOP}, m^{\UOP}, \eta,
\Delta^{\UOP}, \epsilon)$ is a bialgebra in ${}_{H^{op}}^{H^{op}}{\mathcal YD}$. If
$\overline{R}$ is also a bialgebra in ${}_H^H{\mathcal YD}$ and $f : R \longrightarrow
\overline{R}$ is a bialgebra morphism, then $f : R^{\UOP} \longrightarrow
\overline{R}^{\UOP}$ is a bialgebra morphism.

Our assertions about $A^{\UOP}$, $C^{\UOP}$, and $R^{\UOP}$ can be shown directly with
a good deal of effort. For this the $a^{op}$ and $c^{cop}$ notations are strongly
recommended. A more illuminating approach which yields much easier proofs is to
recognize that there is an isomorphism $(F, \vartheta)$ of the monoidal categories
${}_H^H{\mathcal YD}$ and ${}_{H^{op}}^{H^{op}}{\mathcal YD}$. The functor $F :
{}_H^H{\mathcal YD} \longrightarrow {}_{H^{op}}^{H^{op}}{\mathcal YD}$ is defined by
$F(M) = M^{\UOP}$ for objects $M$ and $F(f) = f$ for morphisms $f$. The morphism of
left $H^{op}$-modules $\vartheta_{M, N} : F(M{\otimes}N) \longrightarrow
F(M){\otimes}F(N)$ is defined by $\vartheta_{M, N}(m{\otimes}n) =
S^{-1}(n_{(-1)}){\cdot}m{\otimes}n_{(0)}$ for all $m \in M$ and $n \in N$. Observe
that $\vartheta_{M, N}^{-1} : F(M){\otimes}F(N) \longrightarrow F(M{\otimes}N)$ is
given by $\vartheta_{M, N}^{-1}(m{\otimes}n) = n_{(-1)}{\cdot}m{\otimes}n_{(0)}$ for
all $m \in M$ and $n \in N$.

Let $(A, m, \eta)$ be an algebra in ${}_H^H{\mathcal YD}$. Then $(F(A),
F(m){\circ}\vartheta^{-1}_{A, A}, F(\eta))$ is an algebra of
${}_{H^{op}}^{H^{op}}{\mathcal YD}$. Since $m^{\UOP} = (F(m){\circ}\vartheta_{A,
A}^{-1})^{op}$ it follows that $A^{\UOP}$ is an algebra in
${}_{H^{op}}^{H^{op}}{\mathcal YD}$ as well. Let $(C, \Delta, \epsilon)$ be a
coalgebra of ${}_H^H{\mathcal YD}$. Then $(F(C), \vartheta_{C, C}{\circ}F(\Delta),
F(\epsilon))$ is a coalgebra of ${}_{H^{op}}^{H^{op}}{\mathcal YD}$ which is
$C^{\UOP}$. If $R$ is a bialgebra in ${}_H^H{\mathcal YD}$ then the object $R^{\UOP}$
with its algebra and coalgebra structures $R^{\UOP}$ is a bialgebra in
${}_{H^{op}}^{H^{op}}{\mathcal YD}$.

Suppose that $H$ has bijective antipode and let $V$ by an object of ${}_H^H{\mathcal
YD}$. Using Corollary \ref{CorNicChar} we are able to relate ${\mathfrak B}(V)^{\UOP}$
to a Nichols algebra.

Observe that the grading of ${\mathfrak B}(V)$ is a bialgebra grading for ${\mathfrak
B}(V)^{\UOP}$. It is not hard to see that ${\mathfrak B}(V)(1)$ generates ${\mathfrak
B}(V)^{\UOP}$ and is also the space of primitives of ${\mathfrak B}(V)^{\UOP}$. As a
subobject of ${\mathfrak B}(V)^{\UOP}$ note that ${\mathfrak B}(V)^{\UOP}(1) =
V^{\UOP}$. Thus
\begin{equation}\label{EqBVop}
{\mathfrak B}(V^{\UOP}) = {\mathfrak B}(V)^{\UOP}
\end{equation}
by Corollary \ref{CorNicChar}.
\section{Associated Constructions in ${}_{H^o}^{H^o}{\mathcal YD}$}\label{SecHo}
We now turn to constructions in ${}_{H^o}^{H^o}{\mathcal YD}$. As in the section $H$
has bijective antipode $S$. Starting with objects, algebras, coalgebras, and
bialgebras in ${}_H^H{\mathcal YD}$ we construct counterparts in
${}_{H^o}^{H^o}{\mathcal YD}$ which are important for the analysis of bi-products in
Section \ref{SecUABiProducts}.

First we consider the objects. Let $(M, \eta, \delta) \in {}_H^H{\mathcal YD}$. We
construct an object $(M^r, \delta^o, \eta^o)$ of ${}_{H^o}^{H^o}{\mathcal YD}$. Regard
$H^*{\otimes}M^*$ as a subspace of $(H{\otimes}M)^*$ in the usual way. Recall that
$M^r$, the subspace of all $m^* \in M^*$ which vanish on $I{\cdot}M$ for some cofinite
ideal $I$ of $H$, can be characterized as $M^r = (\eta^*)^{-1}(H^*{\otimes}M^*)$.
Furthermore $\eta^*(M^r) \subseteq H^o{\otimes}M^r$ and  $(M^r, \eta^o)$ is a left
$H^o$-comodule, where $\eta^o = \eta^*|_{M^r}$. Thus the comodule action $\eta^o(m^*)
= m^*_{(-1)}{\otimes}m^*_{(0)}$ of $\eta^o$ on $m^* \in M^r$ is determined by
\begin{equation}\label{EqDeltaO}
m^*_{(-1)}(h)m^*_{(0)}(m) = m^*(h{\cdot}m)
\end{equation}
for all $h \in H$ and $m \in M$. See \cite{CoRef} for example.

The left $H$-comodule structure $(M, \delta)$ induces a (rational) right $H^*$-module
action in $M$ which in turn induces a left $H^*$-module structure $(M^*, {\cdot})$ on
$M^*$ under the transpose action. It is easy to see that $(M^*, {\cdot}) = (M^*,
\delta^*|_{H^*{\otimes}M^*})$. By restriction of the $H^*$-action $M^*$ is a left
$H^o$-module. A straightforward calculation yields
$$
\eta^*(h^o{\cdot}m^*) =
h^o_{\;(1)}m^*_{(-1)}S(h^o_{\;(3)}){\otimes}h^o_{\;(2)}{\cdot}m^*_{(0)}
$$
for all $h^o \in H^o$ and $m^* \in M^r$, where $S$ is the antipode of $H^o$. Let
$\delta^o = \delta^*|_{H^o{\otimes}M^r}$. Thus $H^o{\cdot}M^r \subseteq M^r$; hence
$(M^r, \delta^o)$ is an $H^o$-submodule of $M^*$ which we denote $(M^r, {\cdot})$. The
equivalence of (\ref{EqYDComp}) and (\ref{EqYDCompAntipode}) imply that $(M^r,
{\cdot}, \eta^o) \in {}_{H^o}^{H^o}{\mathcal YD}$. The left $H^o$-module action on
$M^r$ is given explicitly by
\begin{equation}\label{EqDotStar}
(h^o{\cdot}m^*)(m) = m^*(m{\leftharpoonup}h^o) = h^o(m_{(-1)})m^*(m_{(0)})
\end{equation}
for all $h^o \in H^o$, $m^* \in M^r$, and $m \in M$. Note that $(M^{\UOP})^r = M^r$ as
vector spaces. If $N$ is also an object of ${}_H^H{\mathcal YD}$ and $f : M
\longrightarrow N$ is a morphism then $f^*(N^r) \subseteq M^r$ and the restriction
$f^r = f^*|_{N^r}$ is a morphism $f^r : N^r \longrightarrow M^r$ since $f^*(n^*)
(h{\cdot}n) = n^*_{\;\;(-1)}(h)f^*(n^*_{\;\;(0)})(n)$ for all $n^* \in N^r$, $h \in
H$, and $n \in N$.

Suppose that $(C, \Delta, \epsilon)$ is a coalgebra in ${}_H^H{\mathcal YD}$. Then the
object $C^r \in {}_{H^o}^{H^o}{\mathcal YD}$ has the structure of an algebra in the
category; as a $k$-algebra it is a subalgebra of the dual algebra $C^*$. We show that
$C^r$ is a subalgebra of $C^*$ and leave the remaining details of the proof that $C^r$
is an algebra in ${}_{H^o}^{H^o}{\mathcal YD}$ to the reader.

Since $\epsilon : C \longrightarrow k$ is a morphism ${\rm Ker}\,\epsilon$ is a left
$H$-submodule of $C$. Therefore $\epsilon \in C^r$. Suppose that $a, b \in C^r$. Then
$a(I{\cdot}C) = (0) = b(J{\cdot}C)$ for some cofinite ideals $I$ and $J$ of $H$. Since
$\Delta_H$ is an algebra map, $L = \Delta_H^{-1}(I{\otimes}H + H{\otimes}J)$ is an
ideal of $H$ which is cofinite since $I$ and $J$ are cofinite ideals of $H$. Since
$\Delta = \Delta_C$ is a map of left $H$-modules it follows that $\Delta(h{\cdot}c) =
h_{(1)}{\cdot}c^{(1)}{\otimes}h_{(2)}{\cdot}c^{(2)}$ for all $h \in H$ and $c \in C$.
Thus
$$
(ab)(K{\cdot}C) \subseteq (a{\otimes}b)(\Delta(K{\cdot}C)) \subseteq
(a{\otimes}b)(I{\cdot}C{\otimes}C{\cdot}C + C{\cdot}C{\otimes}J{\cdot}C)  = (0)
$$
from which $ab \in C^r$ follows. Also note that if $\overline{C}$ is another coalgebra
in ${}_H^H{\mathcal YD}$ and if $f : C \longrightarrow \overline{C}$ is a coalgebra
morphism then $f^r : \overline{C}^r \longrightarrow C^r$ is an algebra morphism.

One can think of $C^r$ as the counterpart in ${}_{H^o}^{H^o}{\mathcal YD}$ of the dual
k-algebra $C^*$. Suppose $A$ is an algebra in ${}_H^H{\mathcal YD}$. There is a
counterpart $A^{\UO}$ in ${}_{H^o}^{H^o}{\mathcal YD}$ to the dual $k$-coalgebra $A^o$
which arises very naturally in Section \ref{SubAsmashHO}. As a vector space $A^{\UO}$
is the set of all functionals in $A^*$ which vanish on a cofinite subspace $I$ of $A$
which is both an ideal of $A$ and also a left $H$-submodule of $A$.

Now $I^\perp$ is a subcoalgebra of the dual $k$-coalgebra $A^o$ and is also a left
$H$-submodule of $A^r$. Since the intersection of two cofinite subspaces of $A$ which
are both ideals of $A$ and left $H$-submodules of $A$ has the same properties, it
follows that $A^{\UO}$ is a subcoalgebra of $A^o$ and also a subobject of $A^r$. At
this point is not hard to see that $A^{\UO}$ is a coalgebra in
${}_{H^o}^{H^o}{\mathcal YD}$. Let $\overline{A}$ be an algebra in ${}_H^H{\mathcal
YD}$ also and suppose that $f : A \longrightarrow \overline{A}$ is an algebra
morphism. Then $f^r(\overline{A}^{\UO}) \subseteq A^{\UO}$ and the restriction
$f^{\UO} = f^r|_{\overline{A}^{\UO}}$ is a coalgebra morphism $f^{\UO} :
\overline{A}^{\UO} \longrightarrow A^{\UO}$.

Suppose that $R$ is a bialgebra in ${}_H^H{\mathcal YD}$. Then the object $R^{\UO}$ of
${}_{H^o}^{H^o}{\mathcal YD}$ is a bialgebra in ${}_{H^o}^{H^o}{\mathcal YD}$ with the
subalgebra structure of the $k$-algebra $R^*$ and the subcoalgebra structure of the
$k$-coalgebra $R^o$. Furthermore, if $\overline{R}$ is also a bialgebra in
${}_H^H{\mathcal YD}$ and $f : R \longrightarrow \overline{R}$ is a bialgebra
morphism, then $f^{\UO} : \overline{R}^{\UO} \longrightarrow R^{\UO}$ is a bialgebra
morphism.

Let $V$ be an object of ${}_H^H{\mathcal YD}$. There is a natural relationship between
${\mathfrak B}(V)^{\UO}$ and a Nichols algebra. Consider the one-one map $i : V^{r}
\longrightarrow {\mathfrak B}(V)^*$ defined for $v^* \in V^{\UO}$ by
$$
i (v^*) = \left\{\begin{array}{lll} v^*(x) & : & x \in {\mathfrak B}(V)(1) = V \\ 0 &
: & x \in {\mathfrak B}(V)(n), \;\; n \neq 1
\end{array}\right. .
$$
Then ${\rm Im}\, i \subseteq {\mathfrak B}(V)^{\UO}$ and $i : V^{r} \longrightarrow
{\mathfrak B}(V)^{\UO}$ is a one-one morphism. Let ${\mathcal I} : {\mathfrak
B}(V^{r}) \longrightarrow {\mathfrak B}(V)^{\UO}$ be the bialgebra morphism of
Corollary \ref{CorNicChar} which extends $i$. Since ${\rm Ker}\,{\mathcal
I}{\cap}P({\mathfrak B}(V^{r})) = {\rm Ker}\,{\mathcal I}{\cap}V^{r} = (0)$ it
follows that ${\mathcal I}$ is one-one by \cite[Lemma 11.0.1]{SweedlerBook}. again. We
have shown
\begin{equation}\label{EqBVo}
{\mathcal I} : {\mathfrak B}(V^{r}) \longrightarrow {\mathfrak B}(V)^{\UO} \quad
\mbox{is a one-one bialgebra morphism}
\end{equation}
When $V$ is finite-dimensional ${\mathfrak B}(V^{r})$ is identified with the graded
dual of ${\mathfrak B}(V)$ via the map ${\mathcal I}$. In the special case of a Yetter--Drinfel'd module over the group algebra of a finite group with finite-dimensional $\mathfrak{B}(V)$, the dual of $\mathfrak{B}(V)$ was determined in \cite[Theorem 2.2]{B}.
\section{Bilinear Forms in the Yetter--Drinfel'd
Context}\label{SecFormYD}
Let $H$ be a Hopf algebra with bijective antipode. Let $R$, $T$ be bialgebras in
${}_H^H{\mathcal YD}$ and suppose that $\beta : T{\otimes}R \longrightarrow k$ is a
linear form. We will find the following analogs to (A.1)--(A.4) useful:
\begin{enumerate}
\item[{\rm (B.1)}] $\beta (tt', r) = \beta (t,
S^{-1}(r^{(2)}_{\;\;\;(-1)}){\cdot}r^{(1)})\beta (t', r^{(2)}_{\;\;\;(0)})$ for all
$t, t' \in T$ and $r \in R$;
\item[{\rm (B.2)}] $\beta (1, r) = \epsilon(r)$ for all $r \in R$;
\item[{\rm (B.3)}] $\beta (t, rr') = \beta (t^{(2)}, r)\beta
(t^{(1)}, r')$ for all $t \in T$ and $r, r' \in R$; \item[{\rm (B.4)}] $\beta (t, 1) =
\epsilon(t)$ for all $t \in T$.
\end{enumerate}
We leave the reader with the exercise of establishing the equivalence of (B.1)--(B.4)
with analogs of (\ref{EqTauUAopo}) and (\ref{EqTauUAUoop}) :
\begin{Lemma}\label{LemmaBeta}
Let $H$ be a Hopf algebra with bijective antipode over the field $k$, let $T$ and $R$
be bialgebras in ${}_H^H{\mathcal YD}$, and suppose $\beta : T{\otimes}R
\longrightarrow k$ is a linear form. Then the following are equivalent:
\begin{enumerate}
\item[{\rm a)}] (B.1)--(B.4) hold. \item[{\rm b)}] $\beta_\ell (T)
\subseteq R^{\UOP \; \UO}$ and $\beta_\ell : T \longrightarrow R^{\UOP \; \UO}$ is a
bialgebra map. \item[{\rm c)}] $\beta_r (R) \subseteq T^{\UO \; \UOP}$ and $\beta_r :
R \longrightarrow T^{\UO \; \UOP}$ is a bialgebra map.
\end{enumerate}
\epf
\end{Lemma}
\medskip

Let $K$ and $H$ be bialgebras over $k$ and suppose $H$ has bijective antipode $S$. Let
$W \in {}_K^K{\mathcal YD}$, $V \in {}_H^H{\mathcal YD}$, and let $\tau : K{\otimes}H
\longrightarrow k$, $\beta : W{\otimes}V \longrightarrow k$ be linear forms. Two
conditions relating $\tau$ and $\beta$ will play an important part in this paper:
\begin{enumerate}
\item[{\rm (C.1)}] $\beta (k{\cdot}w, v) = \beta(w,
v{\leftharpoonup}\tau_\ell(k))$ for all $k \in K$, $w \in W$, and $v \in V$;
\item[{\rm (C.2)}] $\beta (w{\leftharpoonup}\tau_r(h), v) = \beta(w,
S^{-1}(h){\cdot}v)$ for all $w \in W$, $h \in H$, and $v \in V$.
\end{enumerate}
The first (C.1) implies $V^\perp = {\rm Ker}\,\beta_\ell$ is a $K$-submodule of $W$
and the second (C.2) implies $W^\perp = {\rm Ker}\,\beta_r$ is an $H$-submodule of
$V$. These conditions have formulations in terms of linear and colinear maps.
\begin{Prop}\label{KHWVTauBeta}
Let $K$ and $H$ be Hopf algebras with bijective antipodes over the field $k$. Suppose
that $\tau : K{\otimes}H  \longrightarrow k$ satisfies (A.1)--(A.4), let $W \in
{}_K^K{\mathcal YD}$, $V \in {}_H^H{\mathcal YD}$, and let $\beta : W{\otimes}V
\longrightarrow k$ is a linear form. Then the following are equivalent:
\begin{enumerate}
\item[{\rm a)}] (C.1) and (C.2) hold. \item[{\rm b)}]
$\beta_\ell(W) \subseteq (V^{\UOP})^r$ and $\beta_\ell : W \longrightarrow
(V^{\UOP})^r$ is $\tau_\ell$-linear and colinear.
\item[{\rm c)}] $\beta_r(V) \subseteq W^r$ and $\beta_r : V
\longrightarrow  W^r$ is $\tau_r$-linear and colinear.
\end{enumerate}
\end{Prop}

\pf We show that parts a) and b) are equivalent and leave the reader with the exercise
of adapting our proof to establish the equivalence of parts a) and c). In the latter
the roles of (C.1) and (C.2) are reversed.

Suppose that $\beta_\ell (W) \subseteq (V^{\UOP})^r$ and consider the linear map
$\beta_\ell : W \longrightarrow (V^{\UOP})^r$. Using (\ref{EqDotStar}) it follows that
$\beta_\ell$ is $\tau_\ell$-linear if and only if (C.1) holds. Using (\ref{EqDeltaO}),
where $S^{-1}(h){\cdot}v = h{\cdot}_{op}v$ replaces $h{\cdot}v$,  it follows that
$\beta_\ell$ is $\tau_\ell$-colinear if and only if (C.2) holds.

Suppose that (C.2) holds. Now $W$ is a right $K^*$-module under the rational action
arising from $(M, \delta)$. Now $\tau_r : H \longrightarrow (K^o)^{op}$ is an algebra
map by (\ref{EqTauUAUoop}). Thus $W$ is left $H$-module by pullback along $\tau_\ell$.
Let $w \in W$. Then $H{\cdot}w = w{\leftharpoonup} \tau_r (H)$ is finite-dimensional,
there is a cofinite ideal $I$ of $H$ such that $(0) = I{\cdot}w = w{\leftharpoonup}
\tau_r (I)$. Thus
$$
\beta_\ell(w)(I{\cdot}_{op}V) = \beta(w, S^{-1}(I){\cdot}V) = \beta(w{\leftharpoonup}
\tau_r (I), V) = (0)
$$
which means that $\beta_\ell(w) \in (V^{\UOP})^r$. \epf
\medskip

\begin{Cor}\label{CorBBeta}
Suppose that $K$ and $H$ are Hopf algebras with bijective antipodes over the field $k$
and suppose that $\tau : K{\otimes}H \longrightarrow k$ satisfies (A.1)--(A.4). Let $W
\in {}_K^K{\mathcal YD}$, $W \in {}_H^H{\mathcal YD}$, and $\tau$ and $\beta :
W{\otimes}V \longrightarrow k$ satisfy (C.1) and (C.2). Then:
\begin{enumerate}
\item[{\rm a)}] There is a form ${\mathfrak B}(\beta) : {\mathfrak
B}(W){\otimes}{\mathfrak B}(V) \longrightarrow k$ determined by the properties that it
satisfies (B.1)--(B.4) and ${\mathfrak B}(\beta)|_{W{\otimes}V} = \beta$. Furthermore
${\mathfrak B}(\beta)$ satisfies (C.1) and (C.2). \item[{\rm b)}] Suppose that
$\overline{K}$ and $\overline{H}$ are also Hopf algebras with bijective antipodes
over $k$, $\overline{\tau} : \overline{K}{\otimes}\overline{H} \longrightarrow k$
satisfies (A.1)--(A.4),  $W \in {}_{\overline{K}}^{\overline{K}}{\mathcal YD}$,
$\overline{W} \in {}_{\overline{H}}^{\overline{H}}{\mathcal YD}$, and
$\overline{\tau}$ and $\overline{\beta} : \overline{W}{\otimes}\overline{V}
\longrightarrow k$ satisfy (C.1) and (C.2). If $\overline{\beta}{\circ}(f{\otimes}g) =
\beta$ then ${\mathfrak B}(\overline{\beta}){\circ}({\mathfrak
B}(f){\otimes}{\mathfrak B}(g)) = {\mathfrak B}(\beta)$.
\end{enumerate}
\epf
\end{Cor}
\pf
By Proposition \ref{KHWVTauBeta} b), $\beta_\ell : W \longrightarrow
(V^{\UOP})^r$ is $\tau_\ell$-linear and colinear, and by Corollary \ref{CorNicKH}
$$\mathfrak{B}(\beta_\ell) : \mathfrak{B}(W)\longrightarrow
\mathfrak{B}((V^{\UOP})^r)$$
is a bialgebra map. Then we define ${\mathfrak B}(\beta)_{\ell}$ as the composition of $\mathfrak{B}(\beta_\ell)$ with the maps $$\mathfrak{B}((V^{\UOP})^r)\longrightarrow \mathfrak{B}(V^{\UOP})^{\UO}=\mathfrak{B}(V)^{\UOP \; \UO}$$
in \eqref{EqBVo} and \eqref{EqBVop}. This proves part a), and part b) can be checked easily. 
\epf

\section{Bi-Products Revisited}\label{SectionBiProd}
Let $H$ be a Hopf algebra with antipode $S$ and suppose $R \in {}_H^H{\mathcal YD}$ is
a bialgebra in the category. The biproduct $R{\#}H$ of $R$ and $H$ is a bialgebra over
$k$ described as follows. As a vector space $R{\#}H = R{\otimes}H$ and $r{\#}h$ stands
for the tensor $r{\otimes}h$. As a bialgebra $R{\#}H$ has the smash product and smash
coproduct structures. Thus $1_{R{\#}H} = 1_R{\#}1_H$,
$$
(r{\#}h)(r'{\#}h') = r(h_{(1)}{\cdot}r'){\#}h_{(2)}h'
$$
for all $r, r' \in R$ and $h, h' \in H$,
$$
\Delta (r{\#}h) = (r^{(1)}{\#}r^{(2)}_{\;\;\;
(-1)}h_{(1)}){\otimes}(r^{(2)}_{\;\;\;(0)}{\#}h_{(2)}), \quad \mbox{and} \quad
\epsilon(r{\#}h) = \epsilon(r)\epsilon(h)
$$
for all $r \in R$ and $h \in H$.

The map $j : H \longrightarrow R{\#}H$ defined by $j (h) = 1{\#}h$ for $h \in H$ is a
bialgebra map and the map $\pi : R{\#}H \longrightarrow H$ defined by $\pi (r) =
r{\#}1$ for $r \in R$ is an algebra map which satisfy $\pi{\circ}j = I_H$. Starting
with the bialgebras $A = R{\#}H$, $H$ and the maps $j$, $\pi$ one can recover $R =
R{\#}1$ as a bialgebra in ${}_H^H{\mathcal YD}$. Consider the convolution product
\begin{equation}\label{EqPi}
\Pi = I_A{*}(j{\circ}S{\circ}\pi)
\end{equation}
which as an endomorphism of $A$ given by $\Pi(a) = a_{(1)}j (S(\pi(a_{(2)})))$ for all
$a \in A$. Observe that $\Pi(r{\#}h) = (r{\#}1)\epsilon(h)$ for all $r \in R$ and $h
\in H$. In particular $R = {\rm Im}\,\Pi$. As a $k$-algebra $R$ is merely a subalgebra
of $A$. As a $k$-coalgebra
\begin{equation}\label{EqDeltaSubR}
\Delta_R(r) = \Pi(r_{(1)}){\otimes}r_{(2)} \qquad \mbox{and} \qquad \epsilon_R(r) =
\epsilon(r)
\end{equation}
for all $r \in R$. As an object of ${}_H^H{\mathcal YD}$ the left $H$-module action on
$R$ is given by
\begin{equation}\label{EqRModAct}
h{\cdot}r = j (h_{(1)})rj(S(h_{(2)})
\end{equation}
for all $h \in H$ and $r \in R$ and as a left $H$-comodule action is given by
\begin{equation}\label{EqRComodAct}
\delta(r) = \pi(r_{(1)}){\otimes}r_{(2)}
\end{equation}
for all $r \in R$.

We are in now in a position to look at biproducts in more abstract terms. Let $A$ be a
bialgebra over $k$ and suppose that $j : H \longrightarrow A$, $\pi : A
\longrightarrow H$ are bialgebra maps which satisfy $\pi{\circ}j = I_H$.  Let $\Pi : A
\longrightarrow A$ be defined by (\ref{EqPi}) and set $R = {\rm Im}\,\Pi$. The mapping
$\Pi$ has many important properties which are basic for what follows and which we use
without particular reference in Section \ref{SecUABiProducts}; see \cite{RP} for
example. First of all $\Pi{\circ}\Pi = \Pi$ and $\pi{\circ}\Pi  =
\eta_H{\circ}\epsilon_A$. Since $\Delta(\Pi(a)) =
a_{(1)}(j{\circ}S{\circ}\pi)(a_{(3)}){\otimes}\Pi(a_{(2)})$ for all $a \in A$ it now
follows that $\Delta (R) \subseteq R{\otimes}A$ and
$$
R = {\rm Im}\,\Pi  = \{ a \in A\, | \, a_{(1)}{\otimes}\pi(a_{(2)}) = a{\otimes}1\} =
A^{co \: \pi}.
$$
The last equation is definition. In particular $R$ is a subalgebra of $A$. Since
$\Delta(R) \subseteq A{\otimes}R$ then map $\Delta_R : A \longrightarrow A{\otimes}A$
defined by (\ref{EqDeltaSubR}) satisfies $\Delta_R(R) \subseteq R{\otimes}R$. Using
the fact that $\Pi(aj(h)) = \Pi(a)\epsilon(h)$ for all $a \in A$ and $h \in H$ it
follows by direct calculation that $(R, \Delta_R, \epsilon|_R)$ is a $k$-coalgebra.

Note that $(A, {\cdot}_j, \delta_\pi) \in {}_H^H{\mathcal YD}$, where $h{\cdot}_j a =
j(h_{(1)})aj (S(h_{(2)}))$ for all $h \in H$ and $a \in A$ and $\delta_\pi (a) = \pi
(a_{(1)}){\otimes}a_{(2)}$ for all $a \in A$. Since $\Delta (R) \subseteq A{\otimes}R$
it follows that $R$ is a left $H$-subcomodule of $(A, \delta_\pi)$. Since $h{\cdot}_j
\Pi(a) = \Pi(j (h)a)$ for all $h \in H$ and $a \in A$ we see that $R$ is a left
$H$-submodule of $(A, {\cdot}_j)$. Therefore $R$ is a subobject of $(A, {\cdot}_j,
\delta_\pi)$ and the actions are those described in (\ref{EqRModAct}) and
(\ref{EqRComodAct}). In fact $R$ with these structures is a bialgebra in
${}_H^H{\mathcal YD}$ and the map $R{\#}H \longrightarrow A$ determined by $r{\#}h
\mapsto rj(h)$ is an isomorphism of $k$-bialgebras which we call the canonical
isomorphism. We refer to $R$ with these structures as the bialgebra  in
${}_H^H{\mathcal YD}$ associated to $(A, H, j, \pi)$.

The preceding discussion has been based on a bialgebra $A$ over $k$ with bialgebra
maps $j : H \longrightarrow A$ and $\pi : A \longrightarrow H$ satisfying $\pi{\circ}j
= I_H$. We have the same context for $A^{op}$, $A^{cop}$, thus for $A^{op \; cop}$,
and $A^o$ too. For $j : H^{op} \longrightarrow A^{op}$ and $\pi : A^{op}
\longrightarrow H^{op}$, as well as $j : H^{cop} \longrightarrow A^{cop}$ and $\pi :
A^{cop} \longrightarrow H^{cop}$, are bialgebra maps which satisfy $\pi{\circ}j =
I_H$, and $\pi^o : H^o \longrightarrow A^o$ and $j^o : A^o \longrightarrow H^o$ are
bialgebra maps which satisfy $j^o{\circ}\pi^o = (\pi{\circ}j)^o = I_{H^o}$. It will be
important to us to understand $A^{op}$ and $A^o$ as biproducts. The analysis is rather
detailed and will be carried out in Section \ref{SecUABiProducts}. We do not need to
deal with $A^{cop}$.

We now turn our attention to maps of biproducts. The result we need follows directly
from definitions.
\begin{Prop}\label{PropBiProdMaps}
Let $H$, $\overline{H}$ be Hopf algebras over the field $k$ and let $R \in
{}_H^H{\mathcal YD}$,  $\overline{R} \in {}_{\overline{H}}^{\overline{H}}{\mathcal
YD}$ be bialgebras in their respective categories. Suppose $\varphi : H
\longrightarrow \overline{H}$ is a bialgebra map and $\psi : R \longrightarrow
\overline{R}$ is a map of $k$-algebras and coalgebras which is also $\varphi$-linear
and colinear. Then the linear map $\psi {\#}\varphi : R{\#}H \longrightarrow
\overline{R}{\#}\overline{H}$ defined by $(\psi {\#}\varphi)(r{\#}h) = \psi
(r){\#}\varphi (h)$ for all $r \in R$ and $h \in H$ is a map of bialgebras over $k$.
\epf
\end{Prop}
\medskip

Continuing with the statement of the proposition, note that $(\psi {\#}\varphi){\circ}j =
\overline{j}{\circ}(\psi {\#}\varphi)$ and $\overline{\pi}{\circ}(\psi {\#}\varphi) =
(\psi {\#}\varphi){\circ}\overline{\pi}$. Suppose that $A$, $\overline{A}$ are
bialgebras over $k$. Let  $j : H \longrightarrow A$, $\pi : A \longrightarrow H$ and
$\overline{j} : \overline{H} \longrightarrow \overline{A}$, $\overline{\pi} :
\overline{A} \longrightarrow \overline{H}$ be bialgebra maps which satisfy
$\pi{\circ}j = I_A$, $\overline{\pi}{\circ}\overline{j} = I_{\overline{A}}$
respectively. In light of the proposition a natural requirement for bialgebra maps $f : A
\longrightarrow \overline{A}$ is $\overline{\pi}{\circ}f = f{\circ}\pi$ and $f{\circ}j
= \overline{j}{\circ}f $. When this is the case $f({\rm Im}\,j) \subseteq {\rm
Im}\,\overline{j}$ and $\varphi : H \longrightarrow \overline{H}$ determined by
$\overline{j}{\circ}\varphi = f{\circ}j$ is a bialgebra map, $f(A^{co \; \pi})
\subseteq \overline{A}^{co \; \overline{\pi}}$, and the restriction $f_r = f|_R$ is a
map $f_r : R \longrightarrow \overline{R}$ of $k$-algebras, $k$-coalgebras, and is
$\varphi$-linear and colinear. Furthermore the diagram

\begin{center}
\begin{picture}(100,70)(0,0)
\put(15,0){$A$} \put(5,50){$R{\#}H$} \put(100,0){$\overline{A}$}
\put(90,50){$\overline{R}{\#}\overline{H}$} \put(40,3){\vector(1,0){50}}
\put(45,53){\vector(1,0){35}} \put(20,40){\vector(0,-1){25}}
\put(105,40){\vector(0,-1){25}} \put(60,8){$f$} \put(50,58){$f_r{\#}\varphi$}
\end{picture}
\end{center}
commutes, where the vertical maps are the $k$-bialgebra isomorphisms determined by
$r{\#}h \mapsto rj(h)$ and $\overline{r}{\#}\overline{h} \mapsto
\overline{r}\,\overline{j}(\overline{h})$ respectively.
\section{$(R{\#}H)^{op}$ as a Bi-Product}\label{SubAsmashHOpp}
Throughout this section $H$ is a Hopf algebra with bijective antipode $S$, $A$ is a
bialgebra over $k$, and $j : H \longrightarrow A$, $\pi : A \longrightarrow H$ are
bialgebra maps which satisfy $\pi{\circ}j = I_H$. We will use the results of Section
\ref{SectionBiProd} rather freely and in most cases without particular reference. As
we noted in Section \ref{SecYDCatBiProd} the multiplicative opposite of a bi-product
is a bi-product.

Let $R = A^{co \; \pi}$, let $(R, {\cdot}, \delta)$ be the structure of $R$ as an
object of ${}_H^H{\mathcal YD}$, and let $(R, m, \eta, \Delta, \epsilon)$ the
bialgebra in ${}_H^H{\mathcal YD}$ associated to $(A, H, j, \pi)$. We recall that
$h{\cdot}r = j(h_{(1)})rj(S(h_{(2)}))$ for all $h \in H$ and $r \in R$ by
(\ref{EqRModAct}), and $\delta(r) = \pi(r_{(1)}){\otimes}r_{(2)}$ for all $r \in R$ by
(\ref{EqRComodAct}).

As noted in Section \ref{SectionBiProd} the maps $j : H^{op} \longrightarrow A^{op}$,
$\pi : A^{op} \longrightarrow H^{op}$ are bialgebra maps which satisfy $\pi{\circ}j =
I_{H^{op}}$. We first observe that $R = A^{co \; \pi} = (A^{op})^{co \; \pi}$. Let
$(R, {\cdot}', \delta')$ be the structure of $R = (A^{op})^{co \; \pi}$ as an object
in the category ${}_{H^{op}}^{H^{op}}{\mathcal YD}$ and let $(R, m', \eta', \Delta',
\epsilon')$ be the bialgebra in the category associated with $(A^{op}, H^{op}, j,
\pi)$. The calculation
\begin{eqnarray*}
h^{op}{\cdot}'r
& = & j (h^{op}_{\;\;\;(1)})^{op}r^{op}j(S^{op}(h^{op}_{\;\;\;(2)}))^{op} \\
& = & j (S^{-1}(h_{(2)}))rj (h_{(1)}) \\
& = & j (S^{-1}(h_{(2)}))rj (s(S^{-1}(h_{(1)}))) \\
& = & j (S^{-1}(h)_{(1)})rj (s(S^{-1}(h)_{(2)})) \\
& = & S^{-1}(h){\cdot}r
\end{eqnarray*}
for all $h^{op} = h \in H^{op}$ and $r \in R$ shows that ${\cdot}' = {\cdot}_{op}$.
Since $\delta' =  \delta$ it follows that $(R, {\cdot}', \delta') = (R, {\cdot}_{op},
\delta)$. Thus $R = R^{\UOP}$ as an object of ${}_{H^{op}}^{H^{op}}{\mathcal YD}$.

It is clear that $m' = m^{op}$, $\eta' = \eta$, and $\epsilon' = \epsilon$. To
calculate $\Delta'$ we work from the definition $\Pi^{op} =
I_{A^{op}}*^{op}(j{\circ}S^{op}{\circ}\pi)$ and compute $\Pi^{op}(a) =
((j{\circ}S^{-1}{\circ}\pi)(a_{(2)}))a_{(1)}$ for all $a = a^{op} \in A^{op}$. Now
$\Pi(R) \subseteq A{\otimes}R$ and $\Pi$ acts as the identity on $R$, which thus hold
for $\Pi^{op}$ as well. Since $\Pi(j(h)a) = h{\cdot}_{j}\Pi(a)$ for all $h \in H$ and
$a \in A$, the calculation
\begin{eqnarray*}
\Delta'(r)
& = & \Pi^{op}(r^{op}_{\;\;\;(1)}){\otimes}r^{op}_{\;\;\;(2)} \\
& = & ((j{\circ}S^{-1}{\circ}\pi)(r_{(2)}))r_{(1)}{\otimes}r_{(3)} \\
& = & \Pi(j(S^{-1}(\pi(r_{(2)})))r_{(1)}){\otimes}r_{(3)} \\
& = & S^{-1}(\pi(r_{(2)})){\cdot}_j \Pi(r_{(1)}){\otimes}r_{(3)} \\
& = & S^{-1}(r^{(2)}_{\;\;\;(-1)}){\cdot}r^{(1)}{\otimes}r^{(2)}_{\;\;\;(0)}
\end{eqnarray*}
for all $r \in R$ shows that $\Delta' = \Delta^{\UOP}$. We have shown that the
bialgebra in ${}_{H^{op}}^{H^{op}}{\mathcal YD}$ which is associated to $(A^{op},
H^{op}, j, \pi)$ is $R^{\UOP}$.
\begin{Prop}\label{PropHHDual}
Let $H$ be an Hopf algebra with bijective antipode over $k$, let $A$ be a bialgebra
over $k$, and suppose that $j : H \longrightarrow A$, $\pi : A \longrightarrow H$ are
bialgebra maps which satisfy $\pi{\circ}j = I_H$. Let $R = A^{co \, \pi}$ and let $(R,
m, \eta, \Delta, \epsilon)$ be the bialgebra in ${}_H^H{\mathcal YD}$ associated to
$(A, H, j, \pi)$. Then:
\begin{enumerate}
\item[{\rm a)}] $j : H^{op} \longrightarrow A^{op}$ and $\pi :
A^{op} \longrightarrow H^{op}$ are bialgebra maps which satisfy $\pi{\circ}j =
I_{H^{op}}$, $R^{\UOP} = R$ as a vector space, and the bialgebra in
${}_{H^{op}}^{H^{op}}{\mathcal YD}$ associated to $(A^{op}, H^{op}, j, \pi)$ is
$(R^{\UOP}, m^{\UOP}, \eta, \Delta^{\UOP}, \epsilon)$. \item[{\rm b)}] The map
$\varphi : R^{\UOP} {\#}H^{op} \longrightarrow (R{\#}H)^{op}$ given by
$\varphi(r{\#}h) = (1{\#}h)(r{\#}1)$ for all $r \in R$ and $h \in H$ is an isomorphism
of bialgebras. Furthermore the diagram
\begin{center}
\begin{picture}(100,70)(0,0)
\put(97,53){\vector(0,-1){46}} \put(30,60){\vector(1,0){40}}
\put(24,49.5){\vector(4,-3){60}} \put(-20,57){$R^{\UOP}{\#}H^{op}$}
\put(77,57){$(R{\#}H)^{op}$} \put(90,-5){$A^{op}$} \put(50,65){$\varphi$}
\put(99,30){$f$} \put(50,35){$g$}
\end{picture}
\end{center}
commutes, where $f : R{\#}H \longrightarrow A$ and $g: R^{\UOP}{\#}H^{op}
\longrightarrow A^{op}$ are the canonical isomorphisms.
\end{enumerate}
\end{Prop}

\pf We have established part a). As for part b), we first note that $f(r{\#}h) =
rj(h)$ and $g(r{\#}h) = r^{op}j(h)^{op} = j(h)r$ for all $r \in R$ and $h \in H$.
Therefore $f{\circ}\varphi = g$. This means the diagram commutes and $\varphi =
f^{-1}{\circ}g$ is an isomorphism of bialgebras. \epf
\medskip

Observe that $\varphi^{-1}(r{\#}h) = h_{(1)}{\cdot}_{op}r{\#}h_{(2)} =
S^{-1}(h_{(1)}){\cdot}r{\#}h_{(2)}$ for all $r \in R$ and $h \in H$.

\section{$(R{\#}H)^o$ as a Bi-Product}\label{SubAsmashHO}
As in the preceding section, $H$ is a Hopf algebra with bijective antipode $S$, $A$ is
a bialgebra over $k$, and $j : H \longrightarrow A$, $\pi : A \longrightarrow H$ are
bialgebra maps which satisfy $\pi{\circ}j = I_H$. Again we will use the results of
Section \ref{SectionBiProd} rather freely and in most cases without particular
reference. We observed in Section \ref{SecYDCatBiProd} that the dual of a bi-product
is a bi-product.

As noted in Section \ref{SectionBiProd} the maps $\pi^o : H^o \longrightarrow A^o$,
$j^o : A^o \longrightarrow H^o$ are bialgebra maps which satisfy $j^o{\circ}\pi^o =
I_{H^o}$. We will show that $(A^o)^{co \; j^o}$ can be identified with $R^{\UO}$, find
the structure of $R^{\UO}$ as an object of ${}_{H^o}^{H^o}{\mathcal YD}$,  and then
find its structure as the bialgebra in ${}_{H^o}^{H^o}{\mathcal YD}$ associated to
$(A^o, H^o, \pi^o, j^o)$.

Let $R' = (A^o)^{co \; j^o}$ and $a^o \in A^o$. Then $a^o \in R'$ if and only if
$a^o_{\;\;(1)}{\otimes}j^o(a^o_{\;\;(2)}) = a^o{\otimes}\epsilon$, or equivalently
$a^o(aj(h)) = a^o(a)\epsilon(h)$ for all $a \in A$ and $h \in H$. Recall from Section
\ref{SectionBiProd} that the map $R{\#}H \longrightarrow A$ determined by $r{\#}h
\mapsto rj(h)$ for all $r \in R$ and $h \in H$ is an isomorphism of bialgebras. Since
$A = Rj(H)$ it follows that $a^o \in R'$ if and only if $a^o(rj(h)) =
a^o(r)\epsilon(h)$ for all $r \in R$ and $h \in H$. The isomorphism gives rise to a
linear embedding $i : R^* \longrightarrow A^*$, where $i (r^*)(rj(h)) =
r^*(r)\epsilon(h)$ for all $r^* \in R^*$,  $r \in R$, and $h \in H$. Observe that $R'
\subseteq {\rm Im}\,i$. Thus we can understand $R'$ in terms of $R^*$ via the
embedding.

Our first claim is that $i (R^{\UO}) = R'$. A consequence of the claim is that the
restriction
\begin{equation}\label{EqRUPandRPrime} i|_{R^{\UO}} : R^{\UO}
\longrightarrow R'
\end{equation}
is a linear isomorphism.

To prove our claim, first of all suppose that $r^{\UO} \in R^{\UO}$. To show that $i
(r^{\UO}) \in R'$ we need only show that $i (r^{\UO}) \in A^o$. Since $r^{\UO} \in
R^{\UO}$, by definition $r^{\UO}(J) = (0)$ for some cofinite subspace $J$ of $R$ which
is an ideal and a left $H$-submodule of $R$. We will use the commutation relations
$$
j(h)r = (h_{(1)}{\cdot}_j r)j (h_{(2)}) \qquad \mbox{and} \qquad rj(s(h)) =
j(s(h_{(1)}))(h_{(2)}{\cdot}_j r)
$$
for all $h \in H$ and $r \in R$. Since $S$ is onto and $J$ is a left $H$-subcomodule
of $R$ it follows from the commutations relations that $j(H)J = Jj(H)$. Thus $Jj(H)$
is a left ideal of $A$. Now $i(r^{\UO})$ vanishes on $Jj(H)$ and $Rj(H)^+ =
R(j(H){\cap}{\rm Ker}\,\epsilon)$ as well. Since $H{\cdot}_j R \subseteq R$ and
$j(H)^+$ is a left ideal of $j(H)$, by the first commutation relation $Rj(H)^+$ is a
left ideal of $A$. Since $Jj(H) + Rj(H)^+$ is a cofinite left ideal of $A$ on which $i
(r^{\UO})$ vanishes, it follows that $i (r^{\UO})$ vanishes on a cofinite ideal of
$A$. Thus $i(r^{\UO}) \in A^o$ as required.

Now suppose that $a^o \in R'$. Since $R' \subseteq {\rm Im}\,i$ it follows that $i
(r^*) = a^o$ for some $r^* \in R^*$. By definition $a^o(I) = (0)$ for some cofinite
ideal $I$ of $A$. Since ideals of $A$ are also left $H$-submodules of $A$, and the
subalgebra $R$ of $A$ is also a left $H$-submodule, $J = R{\cap}I$ is a cofinite ideal
of $R$ which is also a left $H$-submodule of $R$. As $(0) = a^o(J) = i (r^*)(J) =
r^*(J)$ we conclude that $r^* \in R^{\UO}$. We have completed the proof of the claim.

Thus $R^{\UO}$ and $R'$ can be identified as vector spaces by the map of
(\ref{EqRUPandRPrime}).  Accordingly we will think of $R'$ as $R^{\UO}$ and show that
$R'$ is $R^{\UO}$ as an object of ${}_{H^o}^{H^o}{\mathcal YD}$ and $R'$ is $R^{\UO}$
as the bialgebra in the category associated to $(A^o, H^o, \pi^o, j^o)$.

Let $h^o \in H^o$ and $r^{\UO} \in R^{\UO}$. The left $H^o$-module structure on $R'$
is given by $h^o{\cdot}_{\pi^o}i(r^{\UO}) =
\pi^o(h^o_{\;\;(1)})i(r^{\UO})\pi^o(S^o(h^o_{\;\;(2)}))$. We evaluate both sides of
this equation at $r \in R$. Since $r_{(1)}{\otimes}\pi(r_{(2)}) = r{\otimes}1$ we have
$\pi(r_{(1)}){\otimes}r_{(2)}{\otimes}\pi(r_{(3)}) =
\pi(r_{(1)}){\otimes}r_{(2)}{\otimes}1$ from which
$\pi(r_{(1)})s(\pi(r_{(3)})){\otimes} r_{(2)} = \pi(r_{(1)}){\otimes}r_{(2)} =
r_{(-1)}{\otimes}r_{(0)}$ follows. Thus
\begin{eqnarray*}
\left(h^o{\cdot}_{\pi^o}i(r^{\UO})\right)(r) & = &
\left(\pi^o(h^o_{\;\;(1)})(r_{(1)})\right)\left(i(r^{\UO})
(r_{(2)})\right)\left(\pi^o(S^o(h^o_{\;\;(2)}))(r_{(3)})\right) \\
& = & \left(h^o_{\;\;(1)}(\pi(r_{(1)}))\right)\left(i(r^{\UO})
(r_{(2)})\right)\left(h^o_{\;\;(2)}(S(\pi(r_{(3)})))\right) \\
& = & h^o\left(\pi(r_{(1)})S(\pi(r_{(3)}))\right)i(r^{\UO})(r_{(2)}) \\
& = & h^o(r_{(-1)})r^{\UO}(r_{(0)}).
\end{eqnarray*}
We have shown that $\left(h^o{\cdot}_{\pi^o}i(r^{\UO})\right)(r) =
h^o(r_{(-1)})r^{\UO}(r_{(0)})$ for all $h^o \in H^o$, $r^{\UO} \in R^{\UO}$, and $r
\in R$.

The left $H^o$-comodule structure on $R'$ is the subcomodule structure afforded by
$(A^o, \delta_{j^o})$. Now $\delta_{j^o}(i(r^{\UO})) = j^o(i(r^{\UO})_{(1)}){\otimes}
i(r^{\UO})_{(2)}$. Using the first commutation relation above we calculate
\begin{eqnarray*}
\delta_{j^o}(i(r^{\UO}))(h{\otimes}r)
& = & i(r^{\UO})_{(1)}(j(h))i(r^{\UO})_{(2)}(r) \\
& = & i(r^{\UO})(j(h)r) \\
& = & i(r^{\UO})((h_{(1)}{\cdot}_{j}r)j(h_{(2)})) \\
& = & r^{\UO}(h_{(1)}{\cdot}_{j}r)\epsilon(h_{(2)}) \\
& = & r^{\UO}(h{\cdot}_{j}r)
\end{eqnarray*}
for all $h \in H$ and $r \in R$. We have shown that $R' = R^{\UO}$ as an object in
${}_{H^o}^{H^o}{\mathcal YD}$.

Next we consider the product in $R'$. Let $r^{\UO}$, $r'^{\UO} \in R^{\UO}$. Since $A
= Rj(H)$ it is easy to see that $i(r^{\UO})$ is determined on $R$. Thus for $r \in R$
the calculation
\begin{eqnarray*}
i(r^{\UO})(r^{(1)})i(r'^{\UO})(r^{(2)})
& = & i(r^{\UO})(\Pi(r_{(1)}))i(r'^{\UO})(r_{(2)}) \\
& = & i(r^{\UO})\left(r_{(1)}j((s{\circ}\pi)(r_{(2)}))\right)i(r'^{\UO})(r_{(3)}) \\
& = & i(r^{\UO})(r_{(1)})\epsilon((s{\circ}\pi)(r_{(2)}))i(r'^{\UO})(r_{(3)}) \\
& = & i(r^{\UO})(r_{(1)})i(r'^{\UO})(r_{(2)}) \\
& = & (i(r^{\UO})i(r'^{\UO}))(r)
\end{eqnarray*}
shows that the product $r^{\UO}r'^{\UO}$ is derived from dual algebra of $(R^{\UO},
\Delta^{\UO})$.

Finally we consider the coproduct of $R'$. First we calculate $\Pi_{A^o}$ in terms of
$\Pi_A = \Pi$. Let $p \in A^o$ and $a \in A$. Then
\begin{eqnarray*}
\left(\Pi_{A^o}(p)\right)(a)
& = & \left(p_{(1)}((\pi^o{\circ}S^o{\circ}j^o)(p_{(2)}))\right)(a) \\
& = & p_{(1)}(a_{(1)})p_{(2)}((j{\circ}s{\circ}\pi)(a_{(2)})) \\
& = & p(a_{(1)}((j{\circ}s{\circ}\pi)(a_{(2)}))) \\
& = & p(\Pi(a))
\end{eqnarray*}
implies $\Pi_{A^o} = (\Pi_A)^o$. Let $r^{\UO} \in R^{\UO}$. By definition
$$\Delta_{R'}(i(r^{\UO})) = \Pi_{A^o}(i(r^{\UO})_{(1)}){\otimes}i(r^{\UO})_{(2)}.
$$
Thus for $r, r' \in R$ we compute
\begin{eqnarray*}
\Delta_{R'}(i(r^{\UO}))(r{\otimes}r')
& = & \left(\Pi_{A^o}(i(r^{\UO})_{(1)})(r)\right)\left(i(r^{\UO})_{(2)}(r')\right) \\
& = & \left(i(r^{\UO})_{(1)}(\Pi(r))\right)\left(i(r^{\UO})_{(2)}(r')\right) \\
& = & i(r^{\UO})(\Pi(r)r') \\
& = & i(r^{\UO})(rr')
\end{eqnarray*}
since $\Pi$ acts as the identity on $R$. Thus the coproduct for $R' = R^{\UO}$ is that
of the dual coalgebra arising from the subalgebra $R$ of $A$. Therefore the bialgebra
$R'$ in the category ${}_{H^o}^{H^o}{\mathcal YD}$ associated to $(A^o, H^o, \pi^o,
j^o)$ is $R^{\UO}$. We regard $R^*{\otimes}H^*$ as a subspace of the vector space
$(R{\#}H)^* = (R{\otimes}H)^*$ in the natural way.
\begin{Prop}\label{PropHHo}
Let $H$ be an Hopf algebra with bijective antipode over $k$, let $A$ be a bialgebra
over $k$ and suppose that $j : H \longrightarrow A$, $\pi : A \longrightarrow H$ are
bialgebra maps which satisfy $\pi{\circ}j = I_H$. Let $R = A^{co \, \pi}$ and let $(R,
m, \eta, \Delta, \epsilon)$ be the bialgebra in ${}_H^H{\mathcal YD}$ associated to
$(A, H, j, \pi)$. Then:
\begin{enumerate}
\item[{\rm a)}] $\pi^o : H^o \longrightarrow A^o$ and $j^o : A^o
\longrightarrow H^o$ are bialgebra maps which satisfy $j^o{\circ}\pi^o = I_{H^{o}}$,
$R^{\UO} = A^{o \; co \, j^o}$ as a vector space under the identification
$r^{\UO}(rj(h)) = r^{\UO}(r)\epsilon(h)$ for all $r \in R$ and $h \in H$, and
$(R^{\UO}, \Delta^{\UO}, \epsilon,  m^{\UO}, \eta)$ is the bialgebra in
${}_{H^o}^{H^o}{\mathcal YD}$ associated to $(A^o, H^o, \pi^o, j^o)$. \item[{\rm b)}]
The map $\vartheta : R^{\UO} {\#}H^o \longrightarrow (R{\#}H)^o$ given by
$\vartheta(r^{\UO}{\#}h^o) = r^{\UO}{\otimes}h^o$ for all $r^{\UO} \in R^{\UO}$ and
$h^o \in H^o$ is an isomorphism of bialgebras. Furthermore the diagram
\begin{center}
\begin{picture}(100,70)(0,0)
\put(0,0){$R'{\#}H^o$} \put(5,50){$R^{\UO}{\#}H^o$} \put(100,0){$A^o$}
\put(85,50){$(R{\#}H)^o$} \put(40,3){\vector(1,0){50}} \put(45,53){\vector(1,0){38}}
\put(20,40){\vector(0,-1){25}} \put(105,40){\vector(0,-1){25}} \put(60,8){$g$}
\put(55,58){$\vartheta$} \put(-30,25){$i|_{R^{\UO}}{\#}I_{H^o}$}
\put(110,25){$(f^{-1})^o$}
\end{picture}
\end{center}
commutes, where $R' = A^{o \; co \, j^o}$, $i|_{R^{\UO}}$ is the linear isomorphism of
(\ref{EqRUPandRPrime}), and the maps $f : R{\#}H \longrightarrow A$, $g : R'{\#}H^o
\longrightarrow A^o$ are the canonical isomorphisms.
\end{enumerate}
\end{Prop}

\pf We have established part a). As for part b), we first note that the structures on
$R^{\UO}$ are due to the identification $i|_{R^{\UO}} : R^{\UO} \longrightarrow R'$.
Thus $i|_{R^{\UO}}{\#}I_{H^o}$ is an isomorphism of bialgebras by Proposition
\ref{PropBiProdMaps}. Now $g$ and $(f^{-1})^o$ are bialgebra isomorphisms also. Thus
part b) will follow once we show that the diagram commutes.

Let $r^{\UO} \in R^{\UO}$, $h^o \in H^o$, $r \in R$, $h \in H$, and set $a = rj(h)$.
Then
\begin{eqnarray*}
\left(\left((f^{-1})^o{\circ}\vartheta\right)(r^{\UO}{\#}h^o)\right)(a)
& = & \vartheta (r^{\UO}{\#}h^o)(f^{-1}(a)) \\
& = & (r^{\UO}{\otimes}h^o)(r{\#}h) \\
& = & r^{\UO}(r)h^o(h).
\end{eqnarray*}
On the other hand
\begin{eqnarray*}
\left(\left(\PH g{\circ}(i|_{R^{\UO}}{\#}I_{H^o})\right)(r^{\UO}{\#}h^o)\right)(a)
& = & \left(\PH g(i (r^{\UO}){\#}h^o)\right)(a) \\
& = & \left(\PH i (r^{\UO})\pi^o(h^o)\right)(rj(h)) \\
& = & i (r^{\UO})(r_{(1)}j(h_{(1)}))\pi^o(h^o)(r_{(2)}j(h_{(2)})) \\
& = & i (r^{\UO})(r_{(1)}j(h_{(1)}))h^o(\pi(r_{(2)}j(h_{(2)}))) \\
& = & r^{\UO}(r_{(1)})\epsilon(h_{(1)})h^o(\pi(r_{(2)})(\pi{\circ}j)(h_{(2)})) \\
& = & r^{\UO}(r_{(1)})h^o(\pi(r_{(2)})h)) \\
& = & r^{\UO}(r)h^o(1h)) \\
& = & r^{\UO}(r)h^o(h).
\end{eqnarray*}
Since $A = Rj(H)$ our calculations show that the diagram of part b) commutes. \epf

\section{Bialgebras ${\mathcal H} = ({\mathcal U}{\otimes}{\mathcal A})^\sigma$
when ${\mathcal U}, {\mathcal A}$ are Bi-Products}\label{SecUABiProducts} Let $K$ and
$H$ be Hopf algebras with bijective antipodes over the field $k$ and suppose that
$U = T{\#}K$, $A = R{\#}H$ are bi-products. We describe an
extensive class of linear forms $(T{\#}K){\otimes}(R{\#}H) \longrightarrow k$ which
satisfy (A.1)--(A.4). All such forms are in bijective correspondence with the
bialgebra maps $T{\#}K \longrightarrow (R{\#}H)^{o \, cop} = (R{\#}H)^{op \, o}$. Our
forms are derived from certain bialgebra maps and determine two-cocycles $\sigma$.

\subsection{The Linear Form $\beta{\#}\tau : (T{\#}K){\otimes}(R{\#}H)
\longrightarrow k$}\label{SubSecForm} We use the isomorphisms of the two preceding
sections to construct bialgebra maps
$$
T{\#}K \longrightarrow (R{\#}H)^{op \; o} \cong R^{\UOP \, \UO}{\#}H^{op \; o}
$$
which determine linear forms $\beta{\#}\tau : (T{\#}K){\otimes}(R{\#}H)
\longrightarrow k$ satisfying (A.1)--(A.4). In the subsequent section we will
investigate the special case when $K$ and $H$ are group algebras of abelian groups.
These are fundamental, interesting in their own right, and arise in representations
theory of pointed Hopf algebras \cite{ASSurvey,DERHJS21}.
\begin{Theorem}\label{ThmTauMain}
Let $K$ and $H$ be Hopf algebras with bijective antipodes over the field $k$. Suppose
that $T$ and $R$ are bialgebras in the categories ${}_K^K{\mathcal YD}$ and
${}_H^H{\mathcal YD}$ respectively and that $\tau : K{\otimes}H \longrightarrow k$ and
$\beta : T{\otimes}R \longrightarrow k$ are linear forms, where $\tau$ satisfies
(A.1)--(A.4), $\beta$ satisfies (B.1)--(B.4), and $\tau$ and $\beta$ satisfy
(C.1)--(C.2). Let $$\beta {\#}\tau : (T{\#}K){\otimes}(R{\#}H) \longrightarrow k$$ be the
linear form determined by
$$
(\beta {\#}\tau) (t{\#}k, r{\#}h) = \beta(t, S^{-1}(h_{(1)}){\cdot}r)\tau(k, h_{(2)})
$$
for all $t \in T$, $k \in K$, $r \in R$, and $h \in H$. Then:
\begin{enumerate}
\item[{\rm a)}] $(\beta {\#}\tau)_\ell$ is the composite
$$
T{\#}K \stackrel{\beta_\ell{\#}\tau_\ell}{\longrightarrow} R^{\UOP \; \UO}{\#}H^{op \;
o} \stackrel{\vartheta}{\longrightarrow} (R^{\UOP}{\#}H^{op})^o
\stackrel{(\varphi^{-1})^o}{\longrightarrow} (R{\#}H)^{op \, o},
$$
where $\vartheta$ and $\varphi$ are defined in part b) of Propositions \ref{PropHHo}
and \ref{PropHHDual} respectively. Thus $\beta {\#}\tau$ satisfies (A.1)--(A.4).
\item[{\rm b)}] $(\beta{\#}\tau)_\ell =
(\varphi^{-1})^*{\circ}(\beta_\ell{\otimes}\tau_\ell)$ and $(\beta{\#}\tau)_r =
(\beta_r{\otimes}\tau_r){\circ}\varphi^{-1}$.
\item[{\rm c)}]
 $\beta {\#}\tau$ is left (respectively right) non-singular if and only
 if $\beta$ and $\tau$ are left (respectively right) non-singular.
\end{enumerate}
\end{Theorem}

\pf Since $\tau$ satisfies (A.1)--(A.4), its equivalent (\ref{EqTauUAopo}), which is
${\rm Im}\,\tau_\ell \subseteq H^{op \, o}$ and $\tau_\ell : K \longrightarrow  H^{op
\, o}$ is a bialgebra map, holds. Likewise, since $\beta$ satisfies (B.1)--(B.4), by
Lemma \ref{LemmaBeta} it follows that ${\rm Im}\,\beta_\ell \subseteq R^{{\UOP}\,
{\UO}}$ and $\beta_\ell : T \longrightarrow R^{{\UOP} \, {\UO}}$ is a map of algebras
and a map of coalgebras. Now $\beta_\ell$ is $\tau_\ell$-linear and colinear since
$\tau$ and $\beta$ satisfy (C.1)--(C.2) by Proposition \ref{KHWVTauBeta}. Therefore
$\beta_\ell {\#}\tau_\ell : T{\#}K \longrightarrow (R^{\UOP})^{\UO}{\#}(H^{op})^o$ is
a bialgebra map by Proposition \ref{PropBiProdMaps}. At this point part a) follows by
Propositions \ref{PropHHDual} and \ref{PropHHo}.

Part b) is a direct consequence of definitions. As for part c) we first note that the
tensor product of two linear maps is one-one if and only if each tensorand is. Since
$\varphi^{-1}$ and $(\varphi^{-1})^*$ are linear isomorphisms, part c) now follows
from part b). \epf
\medskip

Apropos of the theorem, requiring that $\beta{\#}\tau$ be left or right non-singular
seems to be a rather stringent condition. We will find it very natural, and desirable,
for $\beta$ to be non-singular in connection with representation theory. See Section
\ref{SubKHAbelian}.

We shall call a tuple $(K, H, \tau, T, R, \beta)$ which satisfies the hypothesis of
the preceding theorem $2$-{\em cocycle twist datum}. In applications morphisms of
algebras of the type $((T{\#}K){\otimes}(R{\#}H))^\sigma$ will be of interest to us.
As an immediate consequence of Theorem \ref{ThmTauMain}, Propositions
\ref{PropBiProdMaps}, and \cite[Proposition 4.3]{DERHJS21}:
\begin{Cor}\label{Cor2CocycleTwDa}
Let $K$ and $H$ be Hopf algebras with bijective antipodes over the field $k$ and let
$(K, H, \tau, T, R, \beta)$ and $(\overline{K}, \overline{H}, \overline{\tau},
\overline{T}, \overline{R}, \overline{\beta})$ be two-cocycle twist data.  Suppose
that $\varphi : K \longrightarrow \overline{K}$ and $\nu : H \longrightarrow
\overline{H}$ are bialgebra maps, $f : T \longrightarrow \overline{T}$ and $g : R
\longrightarrow \overline{R}$  are algebra and coalgebra maps, where $f$ is
$\varphi$-linear and colinear and $g$ is $\nu$-linear and colinear. Assume further
that $\overline{\tau}{\circ}(\varphi{\otimes}\nu) = \tau$ and
$\overline{\beta}{\circ}(f{\otimes}g) = \beta$.  Then
$(\overline{\beta}{\#}\overline{\tau}){\circ}((f{\#}g){\otimes} (\varphi{\#}\nu)) =
\beta{\#}\tau$. In particular
$$\displaystyle{
(f{\#}\varphi){\otimes}(g{\#}\nu) : \left(\PH (T{\#}K){\otimes}(R{\#}H)\right)^\sigma
\longrightarrow \left(\PH
(\overline{T}{\#}\overline{K}){\otimes}(\overline{R}{\#}\overline{H})\right)^{\overline{\sigma}}
}$$ is a bialgebra map. \qed
\end{Cor}
\subsection{The Fundamental Case $U = {\mathfrak B}(W){\#}K$ and $A = {\mathfrak B}(V){\#}H$.}\label{SecNichols}
We specialize the results of the preceding section to the case of most interest to us:
when $U = {\mathfrak B}(W){\#}K$ and $A = {\mathfrak B}(V){\#}H$ are biproducts of
Nichols algebras with Hopf algebras having bijective antipodes. We are able to express
assumptions involving ${\mathfrak B}(W)$ and ${\mathfrak B}(V)$ in terms of $W$ and
$V$.

\begin{Theorem}\label{ThmNicMain}
Let $K$ and $H$ be Hopf algebras with bijective antipodes over the field $k$ and let
$\tau : K{\otimes}H \longrightarrow k$ be a linear form which satisfies (A.1)--(A.4).
Suppose $W \in {}_K^K{\mathcal YD}$, $V \in {}_H^H{\mathcal YD}$, $\beta :
W{\otimes}V \longrightarrow k$ is a linear form, and $\tau$ and $\beta$ satisfy (C.1)--(C.2). Then:
\begin{enumerate}
\item[{\rm a)}] $(K, H, \tau , {\mathfrak B}(W), {\mathfrak B}(V),
{\mathfrak B}(\beta))$ is a two-cocycle twist datum. \item[{\rm b)}] The linear form
${\mathfrak B}(\beta){\#}\tau : ({\mathfrak B}(W){\#}K){\otimes}({\mathfrak
B}(V){\#}H) \longrightarrow k$ satisfies (A.1)--(A.4).
\end{enumerate}
\end{Theorem}

\pf The form ${\mathfrak B}(\beta) : {\mathfrak B}(W){\otimes}{\mathfrak B}(V)
\longrightarrow k$ of Corollary \ref{CorBBeta} satisfies the hypothesis of Theorem
\ref{ThmTauMain}. \epf
\medskip

We shall call a tuple $(K, H, \tau , W, V, \beta)$ which satisfies the hypothesis of
the preceding theorem {\em Yetter---Drinfel'd} $2$-{\em cocycle twist datum}. By the
preceding theorem if $(K, H, \tau , W, V, \beta)$ is a Yetter--Drinfel'd two-cocycle
twist datum then $(K, H, \tau , {\mathfrak B}(W), {\mathfrak B}(V), {\mathfrak
B}(\beta))$ is a two-cocycle twist datum. We now turn our attention to morphisms.
\begin{Prop}\label{PropNicBiProdMor}
Let $K$, $H$ be Hopf algebras with bijective antipodes over $k$ and $(K, H, \tau , W,
V, \beta)$ and $(\overline{K}, \overline{H}, \overline{\tau} , \overline{W},
\overline{V}, \overline{\beta})$ be Yetter--Drinfel'd two-cocycle twist data. Suppose
that $\varphi : K \longrightarrow  \overline{K}$ and $\nu : H \longrightarrow
\overline{H}$ are bialgebra maps which satisfy
$\overline{\tau}{\circ}(\varphi{\otimes}\nu) = \tau$. Let $f : W \longrightarrow
\overline{W}$ be $\varphi$-linear and colinear, let $g : V \longrightarrow
\overline{V}$ be $\nu$-linear and colinear, and suppose that
$\overline{\beta}{\circ}(f{\otimes}g) = \beta$. Then
$$\displaystyle{
F : \left(\PH ({\mathfrak B}(W){\#}K){\otimes}({\mathfrak B}(V){\#}H)\right)^\sigma
\longrightarrow \left(\PH ({\mathfrak B}(\overline{W}){\#}\overline{K}){\otimes}
({\mathfrak B}(\overline{V}){\#}\overline{H})\right)^{\overline{\sigma}} }
$$
is a bialgebra map, where $F = ({\mathfrak B}(f){\#}\varphi){\otimes}({\mathfrak
B}(g){\#}\nu)$.
\end{Prop}

\pf First of all we observe that $(K, H, \tau , {\mathfrak B}(W), {\mathfrak B}(V),
{\mathfrak B}(\beta))$ and $(\overline{K}, \overline{H}, \overline{\tau} , {\mathfrak
B}(\overline{W}), {\mathfrak B}(\overline{V}), {\mathfrak B}(\overline{\beta}))$ are
two-cocycle twist data by part a) of Theorem \ref{ThmNicMain}. By assumption
$\overline{\tau}{\circ} (\varphi{\otimes}\nu) = \tau$. Since
$\overline{\beta}{\circ}(f{\otimes}g) = \beta$, it follows by Corollary \ref{CorBBeta}
that ${\mathfrak B}(\overline{\beta}){\circ}({\mathfrak B}(f){\otimes}{\mathfrak
B}(g)) = {\mathfrak B}(\beta)$. At this point we apply Corollary \ref{Cor2CocycleTwDa}
to complete the proof. \epf
\medskip

Let $(K, H, \tau , W, V, \beta)$ be a Yetter-Drinfel'd two-cocycle twist datum. Axiom
(C.1) implies that $V^\perp$ is a left $K$-submodule of $W$ and axiom (C.2) implies
that $W^\perp$ is a left $H$-submodule of $V$. Thus the projections $\pi_W: W
\longrightarrow W/V^\perp$ and $\pi_V: V \longrightarrow V/W^\perp$ are module maps.
If $V^\perp$ is a left $K$-subcomodule of $W$ then $\pi_W$ is also a left $K$-comodule
map and likewise if $W^\perp$ is a left $H$-subcomodule of $V$ then $\pi_V$ is also a
left $H$-comodule map.

Suppose that $V^\perp$ and $W^\perp$ are subcomodules. By part c) of Corollary
\ref{CorNicKH} note that ${\mathfrak B}({\pi_W}) : {\mathfrak B}(W) \longrightarrow
{\mathfrak B}(W/V^\perp)$ and ${\mathfrak B}({\pi_W}) : {\mathfrak B}(W)
\longrightarrow {\mathfrak B}(W/V^\perp)$ are onto since $\pi_W$ and $\pi_V$ are.
Recall that the linear form $\overline{\beta} : W/V^\perp{\otimes}V/W^\perp
\longrightarrow k$ determined by $\overline{\beta}{\circ}(\pi_W{\otimes}\pi_V) =
\beta$ is non-singular. With $\varphi$, $\nu$ the identity, $f = \pi_W$, and $g =
\pi_W$, the preceding proposition gives:
\begin{Cor}\label{CorNicNonSing}
Let $K$, $H$ be Hopf algebras with bijective antipodes over $k$ and let $(K, H, \tau ,
W, V, \beta)$ be a Yetter--Drinfel'd two-cocycle twist datum. Suppose that $V^\perp$
and $W^\perp$ are subcomodules of $W$ and $V$ respectively. Let $\pi_W : W
\longrightarrow W/V^\perp$ and $\pi_V : V \longrightarrow V/W^\perp$ be the
projections. Then $(K, H, \tau ,W/V^\perp, V/W^\perp, \overline{\beta})$ is a
Yetter--Drinfel'd two-cocycle twist datum and
$$\displaystyle{
\left(\PH({\mathfrak B}(W){\#}K){\otimes}({\mathfrak B}(V){\#}H)\right)^\sigma
\stackrel{F}{\longrightarrow} \left(\PH ({\mathfrak
B}(W/V^\perp){\#}K){\otimes}({\mathfrak B}(V/W^\perp){\#}H)\right)^{\overline{\sigma}}
}
$$
is a surjective bialgebra map, where $F = ({\mathfrak
B}(\pi_W){\#}I_K){\otimes}({\mathfrak B}(\pi_V){\#}I_H)$. \epf
\end{Cor}
\medskip

The preceding corollary can be used in our study of a class of irreducible
representations of an extensive class of examples of two-cocycle twists in
\cite{DERHJS21}. We will be able to replace the domain of $F$ by its image and thus
assume that $\beta$ is non-singular. The non-singularity of $\beta$ has very
interesting consequences.
\section{The Case when $K$ and $H$ are Group Algebras of Abelian
Groups}\label{SubKHAbelian}
Here we specialize the results of Section \ref{SecNichols} to a very typical case. Let $\Gamma$ be an abelian group. We set
${}_{k[\G]}^{k[\G]}{\mathcal YD} = {\YDG}.$ Suppose that $V \in {\YDG}$, $g \in
\Gamma$, and $\chi \in \widehat{\Gamma}$. We set
$$
V_g = \{v \in V \mid \delta(v) = g{\otimes}v\}
$$
and
$$
V_g^{\chi} = \{ v \in V_g \mid h{\cdot}v = \chi(h)v \; \text{for all}\; h \in
\Gamma \}.
$$
A Yetter-Drinfel'd module in $\YDG$ can be described as a $\Gamma$-graded vector space
which is a $\Gamma$-module such that all $g$-homogeneous components $g \in \Gamma$ are
stable under the $\Gamma$-action.

In this section we fix abelian groups $\Lambda$ and $\Gamma$, positive integers $n$
and $m$, elements $z_1, \ldots, z_n \in \Lambda$ and $g_1,\ldots, g_m \in \Gamma,$ and
nontrivial characters $\eta_1, \ldots, \eta_n \in \widehat{\Lambda}$ and
$\chi_1,\ldots, \chi_m \in \widehat{\Gamma}$. 
Suppose $W \in {\YDL}$ has basis $u_i \in
W_{z_i}^{\eta_i},1 \leq i \leq n$, and that $V \in {\YDG}$ has basis $a_j \in V_{g_j}^{\chi_j},1 \leq j \leq m$. 

The following corollaries will play a useful role in
\cite{DERHJS21}.

\begin{Cor}\label{groupcase}
In addition to the above, let $\varphi : \Lambda \longrightarrow \widehat{\Gamma}$ be a
group homomorphism, $s : \{1, \ldots, n\}
\longrightarrow \{1, \ldots, m\}$ be  a function, and $\lambda_1, \ldots, \lambda_n \in k$.

Let $\tau : k[\Lambda] \otimes k[\G] \longrightarrow k$ and $\beta : W{\otimes}V\longrightarrow k$
be the linear forms defined by 
\begin{equation}\label{Deftau}
\tau(z \otimes g)= \varphi(z)(g) \text{ and }\beta(u_i{\otimes}a_j) = \lambda_i\delta_{s(i), j}
\end{equation}
for all $z \in \Lambda, g \in \Gamma,1 \leq i \leq n$ and $1 \leq j \leq m$. Assume further that for all $1 \leq i \leq n$ with $\lambda_i \neq 0$ and $z \in \Lambda$
\begin{equation}
\varphi(z_i) = \chi^{-1}_{s(i)}\; \mbox{and} \; \eta_i(z) = \varphi(z)(g_{s(i)}).  \label{EqVarPhI2} \\
\end{equation}
Then $(k[\Lambda], k[\G], \tau , W, V, \beta)$ is a Yetter--Drinfel'd two-cocycle twist datum, and the corresponding Hopf algebra map 
$$\Phi = (\mathfrak{B}(\beta) {\#} \tau)_\ell : \mathfrak{B}(W) {\#}k[\Lambda] \longrightarrow (\mathfrak{B}(V){\#}k[\G])^{o\, cop}$$
can be described as follows.

For all $1 \leq i \leq n$ there are a unique algebra map 
\begin{equation}
\gamma_i : \mathfrak{B}(V){\#}k[\G] \longrightarrow k \text{ with }\gamma_i(a_j{\#}1) = 0, \gamma_i(1{\#}g) = \varphi(z_i)(g)
\end{equation}
for all $1 \leq j \leq m$ and $g \in \Gamma,$
and a unique $(\varepsilon,\gamma_i)$-derivation 
\begin{equation}
\delta_i : \mathfrak{B}(V){\#}k[\G] \longrightarrow k \text{ with } \delta_i(a_j{\#}1) = \beta(u_i \otimes a_j), \delta_i(1{\#}g)  = 0
\end{equation}
for all $1 \leq j \leq m$ and $g \in \Gamma.$

Then the algebra map $\Phi : \mathfrak{B}(W) {\#}k[\Lambda] \longrightarrow (\mathfrak{B}(V){\#}k[\G])^{o}$ is determined by
\begin{equation}
\Phi(1{\#}z_i) = \gamma_i, \;\Phi(u_i{\#}1) = \delta_i
\end{equation}
for all $1 \leq i \leq n.$
\end{Cor}
\pf We first show that $(k[\G], k[\Lambda], \tau, W, V,
\beta)$ is Yetter--Drinfel'd two-cocycle datum. Since $H$ is commutative $\tau$
satisfies (\ref{EqTauUAopo}), an equivalent of (A.1)--(A.4). Note that (C.2) holds for
$\tau$ and $\beta$ if and only if for all $1 \leq i \leq n$ and $1 \leq j \leq m$ the
equation $\beta(u_i{\leftharpoonup}\tau_r (g), a_j) = \beta(u_i, S^{-1}(g){\cdot}a_j)$
holds for all $g \in \Gamma$. The latter is $\tau_r(g)(z_i)\beta(u_i, a_j) =
\beta(u_i, \chi_j(g^{-1})a_j)$, or $\tau_\ell(z_i)(g)\lambda_i\delta_{\sigma(i), j} =
\chi_j(g^{-1})\lambda_i\delta_{\sigma(i), j}$, an equivalent of the first equation of
(\ref{EqVarPhI2}). Next we note that (C.1) is equivalent to $\beta(z{\cdot}u_i, a_j) =
\beta(u_i, a_j{\leftharpoonup}\tau_\ell(z))$ which is the same as $\beta(\eta_i(z)u_i,
a_j) = \beta(u_i, \tau_\ell(z)(g_j)a_j)$, or $\eta_i(z)\lambda_i\delta_{\sigma(i), j}
= \varphi(z)(g_j)\lambda_i\delta_{\sigma(i), j}$, for all $z \in \Lambda$, $1 \leq i
\leq n$, and $1 \leq j \leq m$. The latter is equivalent to the second equation of
(\ref{EqVarPhI2}).

We have shown that $(k[\G], k[\Lambda], \tau, W, V, \beta)$ is Yetter--Drinfel'd
two-cocycle datum, and the Corollary   follows by
Theorem \ref{ThmNicMain}.
\epf

\begin{Cor}\label{groupcase2}
Assume the situation of Corollary \ref{groupcase}. Let $I' = \{1 \leq i \leq n \mid \lambda_i \neq 0\},$ and assume that the restriction of $s$ to $I'$ is injective.
Let $V'\subseteq V$ and $W' \subseteq W$ be the Yetter--Drinfel'd submodules with bases $a_{s(i)}, i \in I',$ and $u_i, i \in I'.$ 

Then 
\begin{enumerate}
\item [{\rm a)}] $(k[\Lambda], k[\G], \tau , W', V',
\beta')$ is a Yetter--Drinfel'd two-cocycle twist datum. $V^{\perp} \subseteq W$ and $W^\perp \subseteq V$ are Yetter--Drinfel'd submodules, and the inclusion maps $W' \subseteq W$ and $V' \subseteq V$ define isomorphisms $W' \cong W/V^{\perp}$ and $V' \cong V/W^{\perp}.$  The restriction
$$\beta' : W' \otimes V' \to k$$
of $\beta$ is nondegenerate. 
\item [{\rm b)}] The projections $\pi_W : W \to W', \;\pi_V : V \to V'$ define a surjective bialgebra map
$$
\left(\PH ({\mathfrak B}(W){\#}k[\Lambda]){\otimes}({\mathfrak B}(V){\#}k[\G])\right)^\sigma
\stackrel{F}{\longrightarrow} \left(\PH ({\mathfrak B}(W'){\#}k[\Lambda]){\otimes}({\mathfrak
B}(V'){\#}k[\G])\right)^{\sigma'}
,$$
where
$$
F = ({\mathfrak B}(\pi_W){\#}\id){\otimes}({\mathfrak B}(\pi_V){\#}\id).
$$
\end{enumerate}
\end{Cor}
\pf
By Corollary \ref{groupcase} $(k[\G], k[\Lambda], \tau, W, V, \beta)$ is Yetter--Drinfel'd
two-cocycle datum. Thus $(k[\G], k[\Lambda], \tau, W', V', \beta')$ is
Yetter--Drinfel'd two-cocycle datum as well. 

To show that $\beta'$ is non-singular we compute $W^\perp$ and $V^\perp$. Suppose that $a \in V$ and
write $a = \sum_{j = 1}^mx_j a_j$, where $x_1, \ldots, x_m \in k$. Then $a \in
W^\perp$ if and only if $\beta(u_i, \sum_{j = 1}^mx_j a_j) = 0$, or $\sum_{j = 1}^mx_j
\lambda_i\delta_{s(i), j} = 0$, for all $1 \leq i \leq n$. Therefore $W^\perp$ has
basis $\{a_j \,|\, j \in [m]{\setminus}s(I')\}$. We have shown that $W^\perp {\oplus}V'
= V$.

Let $u \in W$ and write $u = \sum_{i = 1}^n y_i u_i$
where $y_1, \ldots, y_n \in k$. Then $u \in V^\perp$ if and only if $\beta(\sum_{i =
1}^n y_i u_i, a_j) = 0$, or $\sum_{i = 1}^n y_i \lambda_i\delta_{\sigma(i), j} = 0$,
for all $1 \leq j \leq m$. Since $s$ is one-one we conclude that $u \in V^\perp$ if
and only if $y_i\lambda_i = 0$ for all $1 \leq i \leq n$. Thus $V^\perp$ has basis
$\{u_i \,|\, i \in [n]{\setminus}I'\}$. Thus $\beta'$ is non-singular, and part a) is
established.

Note that the maps $\pi_W : W \longrightarrow W'$ and $\pi_V : V \longrightarrow V'$ of
Yetter--Drinfel'd modules can be identified with the projections $W \longrightarrow
W/V^\perp$ and $V \longrightarrow V/W^\perp$ respectively. At this point we apply
Corollary \ref{CorNicNonSing} to complete the proof. \epf

\end{document}